\documentclass[12pt,a4paper,reqno]{amsart}
\usepackage[english]{babel}
\usepackage{amssymb,latexsym,amsfonts,amsthm,upref,amsmath}
\usepackage[margin=1in]{geometry}
\usepackage[foot]{amsaddr}
\usepackage[dvipsnames,x11names]{xcolor}
\definecolor{myblue}{HTML}{003399}
\usepackage{hyperref}
\hypersetup{colorlinks,citecolor=myblue,filecolor=black,linkcolor=myblue,urlcolor=myblue}
\usepackage{comment}
\usepackage{enumerate}
\usepackage{tikz}
\usepackage{float}
\usepackage{cleveref}
\usepackage{mathtools}
\usepackage{cite}
\usepackage{filecontents}
\makeatletter
\newcommand{\leqnomode}{\tagsleft@true}
\newcommand{\reqnomode}{\tagsleft@false}
\newcommand{\cev}[1]{\reflectbox{\ensuremath{\vec{\reflectbox{\ensuremath{#1}}}}}}
\makeatother
\newtheorem*{corintro*}{Corollary}
\newtheorem*{thm*}{Theorem}
\newtheorem*{lem*}{Lemma}
\newtheoremstyle{prim}{}{}{\normalfont}{}{\bfseries}{.}{ }{}
\newtheoremstyle{stil}{}{}{\slshape}{}{\bfseries}{.}{ }{}
\theoremstyle{stil}
\newtheorem{thm}{Theorem}[section]
\newtheoremstyle{defi}{}{}{}{}{\bfseries}{.}{ }{}
\theoremstyle{defi}
\newtheorem{defn}[thm]{Definition}
\theoremstyle{defi}
\newtheorem{rem}[thm]{Remark}
\theoremstyle{stil}
\newtheorem{pro}[thm]{Proposition}
\theoremstyle{stil}
\newtheorem{lem}[thm]{Lemma}
\theoremstyle{stil}
\newtheorem{kor}[thm]{Corollary}
\theoremstyle{prim}

\newenvironment{prf}{\noindent \textit{Proof.}}{\null\hfill$\qed$\hskip
2mm\vskip 2mm}

\newcommand{\DY}{ {\rm DY}}
\newcommand{\DYtld}{ \widetilde{{\rm DY}}}
\newcommand{\DYtldc}{ \widetilde{{\rm DY}}_{\hspace{-1pt}\textrm{crit}}}
\newcommand{\DX}{ {\rm DX}}
\newcommand{\I}{ {\rm I}}

\newcommand{\X}{ {\rm X}}
\newcommand{\V}{\mathcal{V}}

\newcommand{\gr}{ {\rm gr}\ts}

\newcommand{\Y}{ {\rm Y}}

\newcommand{\U}{ {\rm U}}

\newcommand{\R}{ {\overline{R}}}
\newcommand{\Rt}{ {\overline{R'\hspace{-1pt}}}}

\newcommand{\vac}{\hspace{-1pt}\mathop{\mathrm{\boldsymbol{1}}}\hspace{-1pt}}

\newcommand{\g}{\mathfrak{g}}
\newcommand{\on}{\mathfrak{o}}
\newcommand{\spn}{\mathfrak{sp}}

\newcommand{\z}{\mathfrak{z}}

\newcommand{\kp}{\kappa}

\newcommand{\CC}{\mathbb{C}}

\newcommand{\ZZ}{\mathbb{Z}}
\newcommand{\TT}{\mathbb{T}}

\newcommand{\Pc}{\mathcal{P}}

\newcommand{\Sc}{\mathcal{S}}

\newcommand{\Bc}{\mathcal{B}}

\newcommand{\Vc}{\mathcal{V}}\newcommand{\Vcr}{\mathcal{V}_{\textrm{crit}}}

\newcommand{\Dc}{\mathcal{D}}
\newcommand{\ve}{\varepsilon}

\newcommand{\wht}{\widehat}
\newcommand{\wvr}{\overline}

\newcommand{\ot}{\otimes}
\newcommand{\ts}{\hspace{1pt}}

\newcommand{\tr}{ {\rm tr}}

\newcommand{\ndo}{\mathop{\mathrm{End}}}
\newcommand{\om}{\mathop{\mathrm{Hom}}}

\newcommand{\cdotrl}{\mathop{\hspace{-2pt}\underset{\text{RL}}{\cdot}\hspace{-2pt}}}
\newcommand{\cdotlr}{\mathop{\hspace{-2pt}\underset{\text{LR}}{\cdot}\hspace{-2pt}}}

\newcommand{\fand}{\quad\text{and}\quad}
\newcommand{\Fand}{\qquad\text{and}\qquad}

\newcommand{\non}{\nonumber}
\newcommand{\beq}{\begin{equation}}
\newcommand{\eeq}{\end{equation}}
\newcommand{\ben}{\begin{equation*}}
\newcommand{\een}{\end{equation*}}

\begin{document}

\makeatletter
\def\author@andify{%
  \nxandlist {\unskip ,\penalty-1 \space\ignorespaces}%
    {\unskip {} \@@and~}%
    {\unskip \penalty-2 \space \@@and~}%
}
\makeatother

\title{$h$-adic quantum vertex algebras associated with   rational $R$-matrix in types $B$, $C$ and $D$}

\author{Marijana Butorac$^1$}
\address{$^1$ Department of Mathematics, University of Rijeka, Radmile Matej\v{c}i\'{c} 2, 51\,000 Rijeka, Croatia}
\email{mbutorac@math.uniri.hr}

\author{Naihuan Jing$^2$}
\address{$^2$ Department of Mathematics, North Carolina State University, Raleigh, NC 27695, USA}
\email{jing@ncsu.edu}

\author{Slaven Ko\v{z}i\'{c}$^3$}
\address{$^3$ Department of Mathematics, Faculty of Science, University of Zagreb,  Bijeni\v{c}ka cesta 30, 10\,000 Zagreb, Croatia}
\email{kslaven@math.hr}



\begin{abstract}
We introduce the  $h$-adic quantum vertex algebras associated with the rational $R$-matrix in types $B$, $C$ and $D$, thus generalizing the Etingof--Kazhdan's construction in type $A$.
Next, we construct the algebraically independent   generators of the center of the $h$-adic quantum vertex algebra in type $B$ at the critical level, as well as the families of central elements in types $C$ and $D$.
Finally, as an application, we obtain commutative subalgebras of the dual Yangian and the families of central elements of the appropriately completed double Yangian at the critical level, in  types $B$, $C$ and $D$.
\end{abstract}

\maketitle
\section*{Introduction}
\allowdisplaybreaks

The notion of   quantum vertex operator algebra, or, more briefly,  quantum VOA,  was introduced by P. Etingof and D. Kazhdan \cite{EK}. They constructed examples of quantum VOAs associated with the classical $r$-matrix on $\mathfrak{sl}_N$ of rational, trigonometric and elliptic type. The theory of quantum VOAs was further generalized  and  developed by H.-S. Li; see, e.g., \cite{Li0,Li,Li2} and references therein. In particular, Li introduced a certain more general notion of   $h$-adic quantum vertex algebra.  In comparison with quantum VOAs, $h$-adic quantum vertex algebras involve weaker constraints on the braiding operator $\Sc$ which governs the  $\Sc$-locality, a certain quantum version of the locality  property; see \cite{Li} for more details.

In this paper,   following the approach in \cite{EK}, we construct the $h$-adic quantum vertex algebras associated with  rational $R$-matrix in types $B$, $C$ and $D$. As in \cite{EK}, the corresponding vertex operator map is expressed in the form of  quantum currents, which were introduced by N. Yu. Reshetikhin and M. A. Semenov-Tian-Shansky \cite{RS}. Our construction relies on the Poincar\'{e}--Birkhoff--Witt theorem for the double Yangians of the corresponding types, which is due to N. Jing, M. Liu and F. Yang  \cite{JLY}. In particular, the $h$-adic quantum vertex algebra structure is defined on the $h$-adically completed vacuum module over the double Yangian, thus resembling the type $A$ case; cf. \cite{EK,JKMY}.

Next, we study the center of the  $h$-adic quantum vertex algebras in types $B$, $C$ and $D$ at the critical level. In type $B$, we  construct an algebraically independent family of generators of the center.
The classical limit of this family  coincides with the generators of the famous Feigin--Frenkel center \cite{FF}, i.e. of the center of the universal affine vertex algebra in type  $B$  at the critical level, which were found by A. I. Molev \cite{M}.
Furthermore,  we show that the center coincides with the $h$-adically completed polynomial algebra in infinitely many indeterminates, thus resembling the classical case.
In types $C$ and $D$, we also obtain certain algebraically independent families of central elements. However, they do not exhaust the whole center, so they only generate $h$-adically completed polynomial subalgebras of the center. By taking their classical limits we  only reproduce some of the generators of the Feigin--Frenkel center constructed in \cite{M}.
Our construction of central elements relies on a particular case of the fusion procedure for the Brauer algebra, which is due to A. P. Isaev, A. I. Molev and O. V. Ogievetsky \cite{IM,IMO}; see also \cite{Mnew}. It  goes parallel with the corresponding construction \cite{JKMY}, where the center of the Etingof--Kazhdan's quantum VOA in type $A$ was determined, and which relies on the fusion procedure for the symmetric group originated in \cite{J}.

In the end, we show that the aforementioned families of central elements generate commutative subalgebras of the dual Yangians in types $B$, $C$ and $D$,
as suggested by \cite[Remark 11.2.5]{Mnew}. Moreover, by regarding their  images, with respect to the vertex operator map, we find explicit formulae for the
families of central elements of the appropriately completed double Yangians  in  types $B$, $C$ and $D$ at the critical level.

\section{Preliminaries}\label{sec01}
\numberwithin{equation}{section}

\subsection{Affine  Lie algebras in types \texorpdfstring{$B,\,C\text{ and }D$}{B, C and D}}\label{subsec011}

Fix an integer $N\geqslant 2$. Let $\on_N$ be the orthogonal  and let $\spn_N$ be the symplectic Lie algebra, where $N$ is even in the symplectic case.
In order to consider orthogonal and symplectic case simultaneously we denote by $\g_N$ any of the Lie algebras $\on_N$ and  $\spn_N$.
Introduce the scalars $\ve_1,\dots,\ve_N$ by $\ve_1=\ldots =\ve_N =1$ if $\g_N =\on_N$ and by $\ve_1=\ldots =\ve_{N/2} =1$, $\ve_{(N+2)/2}=\ldots =\ve_{N} =-1$ if $\g_N=\spn_N$. For any $i=1,\ldots ,N$ set $i'=N-i+1$. Define the matrix $G=(g_{ij})$ in $\ndo\CC^N$   by $g_{ij}=\delta_{ij'}\ve_i$ for all $i,j=1,\ldots ,N$. Clearly,
$G$ is symmetric in the  orthogonal and   skew-symmetric in the symplectic case.
For any matrix $A=(a_{ij})$ in $\ndo\CC^N$ let $A'=G A^t G^{-1}$, where $A^t$ denotes the transposed matrix $A^t=(a_{ji})$. We have $A'=(\ve_i \ve_j a_{j'i'})$.
Set
$$
\sigma =\begin{cases}
1,&\text{if }\g_N=\on_N,\\
2,&\text{if }\g_N=\spn_N.
\end{cases}
$$

Define the operators $P$ and $Q$ in $\ndo\CC^N \ot\ndo\CC^N$   by
$$P=\sum_{i,j=1}^N e_{ij}\ot e_{ji}\fand Q=\sum_{i,j=1}^N \ve_{i}\ve_j \ts e_{ij}\ot e_{i'j'},$$
where $e_{ij}$ are matrix units.
One can easily verify that
\beq\label{formulae2}
P^{t_1} =  P^{t_2},\quad  P'^{{}_1} =   P'^{{}_2},\quad Q^{t_1} =  Q^{t_2},\quad Q'^{{}_1} =   Q'^{{}_2},
\eeq
where    ${}^{t_i}$ and ${}'^{{}_i}$ denote the transpositions ${}^t$ and  ${}'^{{}}$ applied on the $i$-th  factor of the tensor product algebra $\ndo\CC^N \ot \ndo\CC^N$ for $i=1,2$. Therefore, we will usually omit the index $i$ and denote the expressions in \eqref{formulae2} by $P^t, P', Q^t,Q'$ respectively.
The operators $P$ and $Q$ posses the following properties:
\begin{align}\label{formulae}
P'=Q,\quad Q'=P,\quad
P^2 =1,\quad Q^2 =NQ,\quad PQ=QP=\begin{cases}
Q, &\text{if }\g_N =\on_N,\\
-Q,&\text{if }\g_N=\spn_N .
\end{cases}
\end{align}

Recall that the universal enveloping algebra $\U(\g_N)$ is the associative algebra generated by the elements $f_{ij}$, where $i,j=1,\ldots ,N$, subject to the defining relations
\begin{align}\label{unifin}
&F_1 F_2 -F_2 F_1 =(P-Q)F_2 -F_2 (P-Q)\fand F+F'=0.
\end{align}
The elements $F$ and $F'$ in $\ndo\CC^N \ot \U(\g_N)$ are given by
$$F=\sum_{i,j=1}^N e_{ij}\ot f_{ij}\Fand F'=\sum_{i,j=1}^N  e_{ij}\ot \ve_i\ve_j \ts f_{j'i'}.$$
The copies of the matrix in the tensor product algebra $(\ndo\CC^N)^{\ot 2}\ot \U(\g_N)$ are indicated by the subscripts, so that in \eqref{unifin} we have
$$F_1=\sum_{i,j=1}^N e_{ij}\ot1\ot f_{ij}\Fand F_2=\sum_{i,j=1}^N 1\ot e_{ij}\ot f_{ij}.$$

Consider the affine Kac--Moody Lie algebra $\wht{\g}_N=\g_N\ot\CC[t,t^{-1}]\oplus\CC C$; see \cite{Kac} for more details.
Its universal enveloping algebra $\U(\wht{\g}_N)$ is generated by the central element $C$ and the elements $f_{ij}(r)=f_{ij}\ot t^r$, where $i,j=1,\ldots ,N$ and $r\in\mathbb{Z}$, subject to the defining relations
\begin{alignat}{3}
&F(r)_1 F(s)_2 -F(s)_2 F(r)_1=&&\,(P-Q)F(r+s)_2 -F(r+s)_2(P-Q) \non\\
& &&\, +\sigma r \delta_{r+s\ts 0} (P-Q)C,\label{uniaf}\\
& F(r) + F(r)'=0 .&&
\end{alignat}
The elements $F(r)$ and $F(r)'$ in $\ndo\CC^N \ot \U(\wht{\g}_N)$ are defined by
$$F(r)=\sum_{i,j=1}^N e_{ij}\ot f_{ij}(r)\Fand F(r)'=\sum_{i,j=1}^N  e_{ij}\ot \ve_i\ve_j \ts f_{j'i'}(r).$$
As in \eqref{unifin}, the subscripts in \eqref{uniaf}  indicate the copy of $\ndo\CC^N$ in the tensor product algebra $(\ndo\CC^N)^{\ot 2} \ot \U(\wht{\g}_N)$.
Defining relations \eqref{uniaf} can be equivalently written as
\begin{align}
&F_1(u) F_2(v) -F_2(v) F_1(u)  =
\left((P-Q) F_2(v)-F_2(v)(P-Q)\right)\frac{1}{v}\delta\left(\frac{u}{v}\right) \non\\
&\qquad\qquad\qquad\qquad\qquad\qquad-\sigma (P-Q) \ts C\ts\frac{1}{v}\frac{\partial}{\partial u}\delta\left(\frac{u}{v}\right),
\label{uniaf3}
\end{align}
where $\delta(z)=\sum_{r\in\ZZ} z^r \in\CC[[z^{\pm 1}]]$ is the formal delta function
and
$$F(u)=\sum_{r\in\ZZ} F(-r)u^{r-1}\in \ndo\CC^N \ot \U(\wht{\g}_N )[[u^{\pm 1}]] .$$

Introduce the series
\begin{align*}
&F^+(u)=\sum_{r\geqslant 1} F(-r)u^{r-1} \in \ndo\CC^N \ot \U(t^{-1}\g_N [t^{-1}])[[u]],\\
&F^-(u)=\sum_{r\leqslant 0} F(-r)u^{r-1} \in  \ndo\CC^N \ot u^{-1}\U(\g_N [t])[[u^{-1}]],
\end{align*}
so that $F(u)=F^+ (u) +F^-(u)$. Relation \eqref{uniaf3} implies
\begin{align}
F_1^-(u) F_2^+(v) -F_2^+ (v) F_1^-(u) =&\,\frac{1}{u-v} \Big(
(P-Q)\left( F_2^+(v)-F_1^-(u)\right)  \Big.  \label{uniaf4}\\
&\Big.  \,-\left( F_2^+(v)-F_1^-(u)\right)(P-Q)
\Big) + \frac{\sigma C}{(u-v)^{2}} \ts (P-Q) . \non
\end{align}
Throughout this paper, we employ the following expansion convention. For any variables $u_1,\ldots , u_k$ expressions of the form $(u_1+\ldots +u_k)^r$ with $r<0$ should be expanded in nonnegative powers of the variables $u_2,\ldots ,u_k$. For example, in \eqref{uniaf4} we have $k=2$, $r=-2,-1$ and
$$(u-v)^{r}=\sum_{k\geqslant 0} \binom{r}{k} u^{r-k} v^k.$$

Recall that the vacuum module $V_c(\g_N)$ of level $c\in\CC$ for the  algebra $\U(\wht{\g}_N)$ is isomorphic to the
universal enveloping algebra $\U(t^{-1}\g_N [t^{-1}])$   as a complex vector space.
The algebra $\U(t^{-1}\g_N [t^{-1}])$ is generated by the elements $f_{ij}(r)=f_{ij}\ot t^r$, where $i,j=1,\ldots ,N$ and $r\in\mathbb{Z}_{<0}$, subject to the defining relations
\begin{alignat}{3}
&F_1^+(u) F_2^+(v) -F_2^+(v) F_1^+(u) =&&\, \frac{1}{u-v} \Big(\left(F_1^+(u)+ F_2^+(v)\right)(P-Q)\Big.\non\\
& &&\Big.\, -(P-Q)\left(F_1^+(u)+ F_2^+(v)\right)\Big)
\label{vacuum3}
,\\&F^+(u) + F^+(u)'=0
.
\label{vacuum4}
\end{alignat}

\subsection{Rational $R$-matrix}\label{subsec021}

Let $h$ be a formal parameter.
Consider the rational {\em $R$-matrix} $R(u)$, as defined in \cite{ZZ}, over the commutative  ring $\CC[[h]]$,
\beq\label{R}
R(u)=1 -\frac{hP}{u} + \frac{hQ}{u-h\kp},
\eeq
where
$$
 \kp=\begin{cases}
N/2-1,&\text{if }\g_N=\on_N,\\
N/2+1,&\text{if }\g_N=\spn_N.
\end{cases}
$$
$R$-matrix \eqref{R} satisfies the Yang--Baxter equation
\beq\label{YBE}
R_{12}(u)R_{13}(u+v)R_{23}(v)=R_{23}(v)R_{13}(u+v)R_{12}(u).
\eeq
Both sides of \eqref{YBE} are regarded as operators on the triple tensor product $(\CC^N)^{\ot 3}$ and the subscripts
indicate the copies of $\CC^N$ on which the $R$-matrices are applied. For example, we have $R_{23}(u)= 1\ot R(u) $.

By \eqref{formulae2} we have $R(u)^{t_1}=R(u)^{t_2}$ and  $R(u)'^{{}_1}=R(u)'^{{}_2}$, so we   denote these transposed $R$-matrices by  $R(u)^{t}$ and $R(u)'$ respectively.
Using properties \eqref{formulae}
one can easily prove
\beq\label{almost}
  R(u)R(u+h\kp)'=1-h^2 u^{-2}\fand R(u)R(-u)=1- h^2 u^{-2}.
\eeq

\begin{lem}
There exists a unique series $f(u)\in 1+u^{-2}\CC[[u^{-1}]]$ satisfying
\beq\label{normalize}
f(u)f(u+\kp)=\left(1-u^{-2}\right)^{-1}\Fand f(u)f(-u)=\left(1-u^{-2}\right)^{-1}.
\eeq
\end{lem}

\begin{prf}
Write the series $f$ as $f(u)=1+\sum_{r\geqslant 1} f_r u^{-r}$ for some scalars $f_r\in\CC$. The first equality in \eqref{normalize} implies
\beq\label{norm3}
\left(1+\sum_{r\geqslant 1} f_r u^{-r}\right)
\left(1+\sum_{r\geqslant 1}  \left(\sum_{s=1}^r \binom{-s}{r-s} \kp^{r-s} f_s\right) u^{-r}\right)
=\sum_{r\geqslant 0} u^{-2r}.
\eeq
One can easily see that the coefficients $f_r$ are uniquely determined by \eqref{norm3} and that the first few terms of the series $f(u)$ are found by
\beq\label{formoff}
f(u) = 1 + \frac{1}{2}u^{-2} + \frac{\kp}{2}u^{-3} +\frac{3}{8} u^{-4} +\ldots
\eeq

It remains to prove that the series $f(u)=1+\sum_{r\geqslant 1} f_r u^{-r}$ satisfies the second equality in \eqref{normalize}.
Consider the series
$$F(u)=f(u)f(-u)(1-u^{-2})\,\in\, 1+u^{-2}\CC[[u^{-1}]].$$
Multiplying the first equality in \eqref{normalize} by
$f(-u-\kp)f(-u)=\left(1-\left(u+\kp\right)^{-2}\right)^{-1}$  we find
$F(u)=F(u+\kp)^{-1}$. This clearly  implies $F(u)=F(u+2\kp)$.
However, the only series  in $1+u^{-2}\CC[[u^{-1}]]$  which satisfies that equality is the constant series $1$, so we conclude that $f(u)f(-u)(1-u^{-2})=1$, as required.
\end{prf}

As with the $R$-matrix $R(u)$,  the normalized $R$-matrix $\R(u)=f(u/h) R(u)$  satisfies   Yang--Baxter equation \eqref{YBE}.
Furthermore, by the first equalities in \eqref{almost} and \eqref{normalize}   it possesses the {\em crossing symmetry} properties
\beq\label{csym}
\R(u)\ts \R(u+h\kp  )' =1\fand  \R(u+h\kp  )'\ts \R(u) =1.
\eeq
Using the ordered product notation we rewrite \eqref{csym} as
\beq\label{csym2}
\R(u)'\cdotrl \R(u+h\kp  ) =1\fand  \R(u)'\cdotlr \R(u+h\kp  ) =1,
\eeq
where the subscript RL (LR) in \eqref{csym2} indicates that the first tensor factor of $\R(u)'$ is applied from the right (left) while the second tensor factor  of $\R(u)'$ is applied from the left (right). Note that the equations in \eqref{csym2} can be equivalently written as
$$\left(\left(\R(u)'\right)^t\ts \R(u+h\kp  )^t\right)^{t_1} =1\fand  \left(\left(\R(u)'\right)^t\ts \R(u+h\kp  )^t\right)^{t_2} =1.$$
By the second equalities in \eqref{almost} and \eqref{normalize}  the $R$-matrix $\R(u) $ possesses the {\em unitarity} property
\beq\label{uni}
\R(u)\ts \R(-u)=\R(-u)\ts \R(u) =1.
\eeq
Hence we also have
\beq\label{csym3}
\R(-u)' = \R(u+h\kp).
\eeq

\subsection{Double Yangians in types \texorpdfstring{$B,\,C\text{ and }D$}{B, C and D}}\label{subsec022}

We  follow \cite{JLY} to introduce the double  Yangian   in types $B$, $C$ and $D$. In contrast with \cite{JLY}, where the (extended) Yangian double is defined as an associative algebra over $\CC$, all algebras in this paper are defined over the commutative ring $\CC[[h]]$. Such modification  makes the structure of double Yangian  suitable for the construction of $h$-adic quantum vertex algebras.

The  {\em extended Yangian double} $\DX(\g_N)$ for $\g_N$  is defined as the associative algebra over the ring $\CC[[h]]$ generated by the central element $C$ and   elements $t_{ij}^{(r)}$ and $t_{ij}^{(-r)}$, where $i,j=1,\ldots ,N$ and $r=1,2,\ldots$, subject to the defining relations
\begin{align}
R(u-v)T_1(u)T_2(v)&=T_2(v)T_1(u)R(u-v),\label{DX1}\\
R(u-v)T_1^+(u)T_2^+(v)&=T_2^+(v)T_1^+(u)R(u-v),\label{DX2}\\
\R(u-v+h\sigma C/2)T_1(u)T_2^+(v)&=T_2^+(v)T_1(u)\R(u-v-h\sigma C/2).\label{DX3}
\end{align}
The elements $T(u)$ and $T^+(u)$ in $\ndo\CC^N \ot \DX(\g_N)[[u^{\pm 1}]]$ are defined by
\begin{align*}
T(u)=\sum_{i,j=1}^N e_{ij}\ot t_{ij}(u)\Fand T^{+}(u)=\sum_{i,j=1}^N e_{ij}\ot t_{ij}^{+}(u),
\end{align*}
where the $e_{ij}$ denote the matrix units
and the series $t_{ij}(u) $ and $t_{ij}^+ (u) $ are given by
$$t_{ij}(u)=\delta_{ij}+h\sum_{r=1}^{\infty} t_{ij}^{(r)}u^{-r}\Fand t_{ij}^+ (u)=\delta_{ij}-h\sum_{r=1}^{\infty} t_{ij}^{(-r)}u^{r-1}.$$

The {\em double Yangian} $\DY(\g_N)$ for $\g_N$ is defined as the quotient of the $h$-adically completed algebra $\DX(\g_N)[[h]]$ by the relations
\beq\label{DY1}
 T(u) T(u+h\kp)'=1
\fand T^+(u) T^+(u+h\kp)'=1.
\eeq

We now introduce certain  subalgebras of the double Yangian; cf. \cite[Chapter 11]{Mnew}. The {\em  Yangian} $\Y(\g_N)$ is defined as the subalgebra of $\DY(\g_N)$ generated by   all elements $t_{ij}^{(r)}$, where $i,j=1,\ldots ,N$ and $r\in\ZZ$.  By the Poincar\'{e}--Birkhoff--Witt theorem for the double Yangian \cite{JLY},   $\Y(\g_N)$ is isomorphic to the associative algebra over $\CC[[h]]$ generated by  the elements $t_{ij}^{(r)}$, where $i,j=1,\ldots ,N$ and $r\in\ZZ$, subject to  defining relations \eqref{DX1} and
$$T(u) T(u+h\kp)'=1.$$

The {\em dual  Yangian} $\Y^+(\g_N)$ is defined as the  $h$-adically completed subalgebra of $\DY(\g_N)$ generated by   all elements $t_{ij}^{(-r)}$, where $i,j=1,\ldots ,N$ and $r\in\ZZ$.
Define the {\em extended dual Yangian} $\X^+(\g_N)$ as the associative algebra over $\CC[[h]]$ generated by the elements
$t_{ij}^{(-r)}$, where $i,j=1,\ldots ,N$ and $r=1,2,\ldots$, subject to defining relations \eqref{DX2}.
By the Poincar\'{e}--Birkhoff--Witt theorem for the double Yangian \cite{JLY},   $\Y^+(\g_N)$   is isomorphic to the quotient of the $h$-adically completed extended dual Yangian $\X^+(\g_N)[[h]]$
 by the relation
\beq\label{TT}
  T^+(u) T^+(u+h\kp)'=1.
\eeq

For any $c\in\CC$ define the {\em double Yangian $\DY_c(\g_N)$ at the level} $c$ as the quotient of $\DY(\g_N)$ by the ideal generated by $C-c$. The {\em vacuum module $\Vc_c(\g_N)$ at the level} $c$ is defined as the $h$-adic completion of the quotient
$
\DY_c(\g_N) / \I
$, where $\I$ denotes the $h$-adically completed left ideal in $\DY_c(\g_N)$ generated by the elements $t_{ij}^{(r)}$ with $i,j=1,\ldots ,N$ and $r\geqslant 1$. By the Poincar\'{e}--Birkhoff--Witt theorem for the double Yangian \cite{JLY},  $\Vc_c(\g_N)$ and $\Y^+(\g_N)$ are isomorphic as  $\CC[[h]]$-modules. Moreover, one can easily verify that
 the vacuum module $\Vc_c(\g_N)$ is topologically free.

As in \cite{JLY} and \cite[Section 11.2]{Mnew}, set
\beq\label{grading1}
\deg 1=0\Fand  \deg t_{ij}^{(-r)}=-r\quad\text{for all }i,j=1,\ldots ,N,\, r\geqslant 1
\eeq
 and denote by $\Y^{(-r)}$ the  $\CC[[h]]$-span of the elements of $\Y^+(\g_N)$ whose degrees do not exceed $-r$.
We have the ascending filtration
\beq\label{filtration}
  \ldots\subset\Y^{(-r-1)}\subset\Y^{(-r)}\subset\ldots\subset\Y^{(-2)}\subset\Y^{(-1)}\subset\Y^{(0)}=\Y^+(\g_N)  .
	\eeq
Consider the associated graded algebra
$$\gr  \Y^+(\g_N)=\bigoplus_{r\geqslant 0} \Y^{(-r)}/ \Y^{(-r-1)}.$$
Observe that the $\CC[[h]]$-algebra $\gr  \Y^+(\g_N)$ is no longer $h$-adically complete.

By employing the notation
$$t(u)=h^{-1}\left(T(u)-1\right)\fand
t^+(u)=h^{-1} \left(1-T^+ (u)\right)$$
 we express  defining relation \eqref{DX2} for the dual Yangian $\Y^+(\g_N)$ as
\begin{align}
t_1^+ (u)\ts t_2^+ (v) -t_2^+ (v)\ts t_1^+ (u)
= &\ts\ts\frac{1}{u-v}\left(\left(t_1^+ (u) + t_2^+ (v)\right)P -P\left(t_1^+ (u) + t_2^+ (v)\right) \right)\non\\
&-\frac{1}{u-v-h\kp}\left(\left(t_1^+ (u) + t_2^+ (v)\right)Q -Q\left(t_1^+ (u) + t_2^+ (v)\right) \right)\non\\
&+\frac{h}{u-v}\left(P\ts t_1^+ (u) \ts t_2^+ (v) -   t_2^+ (v)\ts t_1^+ (u)\ts P \right)\non\\
&-\frac{h}{u-v-h\kp}\left(Q\ts t_1^+ (u)\ts  t_2^+ (v) -   t_2^+ (v)\ts t_1^+ (u)\ts Q\right).\label{dualyangian}
\end{align}
 Furthermore, defining relation \eqref{TT} can be written as
\begin{align}
t^+(u) + t^+(u+h\kp)'=h\ts t^+(u) \ts t^+(u+h\kp)'. \label{dualyangian2}
\end{align}
Clearly, the relations obtained by taking the highest degree terms in   \eqref{dualyangian} and \eqref{dualyangian2}, with respect to degree operator \eqref{grading1}, coincide with relations
\eqref{vacuum3} and \eqref{vacuum4} respectively. Moreover, from defining relations \eqref{DX3}  we derive the following formula for the action of the Yangian generators  on   $t^+(v)\vac\in\Vc_c (\g_N)[[v]]$:
\begin{align}
t_{1}(u) t_{2}^+ (v)\vac=\,
&\frac{1}{u-v}
\left((P-Q)t_2^+(u)\vac-t_2^+(u)\vac(P-Q) \right)\non\\
&+ \frac{\sigma c}{(u-v)^2}\ts (P-Q) +\,\ldots,\label{dualyangian3}
\end{align}
where the ellipsis  represents the summands of lower degree  with respect to \eqref{grading1}. The relations obtained by taking the highest degree terms in  \eqref{dualyangian3}, with respect to  \eqref{grading1}, coincide with relations obtained by applying
\eqref{uniaf4} to the vacuum vector $\vac\in V_c(\g_N)$.

Denote by $\wvr{t}_{ij}^{(-r)}$ the image of the element $t_{ij}^{(-r)}$ in the $(-r)$-th component of $\gr  \Y^+(\g_N)$.
 The next proposition  is required in the proof of Theorem \ref{FFthm}.  It is a consequence of the Poincar\'{e}--Birkhoff--Witt theorem for the double Yangians of type $B$, $C$ and $D$; see \cite{JLY}.

\begin{pro}\label{pbwthm}
\emph{(a)}
The assignments
\beq\label{pbw}
\wvr{t}_{ij}^{(-r)}\ts\mapsto\ts f_{ij}\ot t^{-r}
\eeq
with $r\geqslant 1$ and $i,j=1,\ldots ,N$
define an isomorphism of $\CC[[h]]$-algebras
$$\gr  \Y^+(\g_N)\ts\cong\ts \U(t^{-1}\g_N [t^{-1}])\ot_{\CC} \CC[[h]].$$
\emph{(b)} For any $c\in\CC$ and integer $n\geqslant 1$ the image of
\beq\label{ontheelements}
t_{0} (u)\left(  t_{1}^+ (v_1)\ldots t_{n}^{+}(v_n)\right) \,\in\,(\ndo\CC^N)^{\ot (n+1)}\ot \Vc_c(\g_N) [[u^{-1},v_1,\ldots,v_n]]
\eeq
with respect to map \eqref{pbw} is equal to
$$ F_{0}^{-}(u) \left(F_{1}^+ (v_1)\ldots F_{n}^{+}(v_n)\right) \in (\ndo\CC^N)^{\ot (n+1)}\ot  \U(t^{-1}\g_N[t^{-1}])[[ u^{-1},v_1,\ldots,v_n]].$$
\end{pro}

\begin{rem}
It is worth noting that defining relation \eqref{DX3} differs from the corresponding relation in \cite{JLY} because we use the normalized $R$-matrix $\R (u)$ instead of $R(u)$. However, this does not affect the action of map \eqref{pbw} on   elements \eqref{ontheelements}. More specifically, due to the form of the normalizing function \eqref{formoff}, this produces only additional terms of the lower degree  with respect to   \eqref{grading1}, which are annihilated by map \eqref{pbw}.
\end{rem}

Let $m\geqslant 1$ be an arbitrary integer, $v=(v_1,\ldots, v_m)$ an $m$-tuple of variables and $z$ a single variable.
Label the tensor factors of $(\ndo\mathbb{C}^{N})^{\ot  (n+m)}$
as follows,
\beq\label{nhk78}\overbrace{(\ndo\mathbb{C}^{N})^{\ot  n}}^{1} \otimes
\overbrace{(\ndo\mathbb{C}^{N})^{\ot  m}}^{2}.
\eeq
Introduce the  functions  with values in
\eqref{nhk78}
by
\begin{align}
\R_{nm}^{12}(u|v|z)= \prod_{i=1,\dots,n}^{\longrightarrow}
\prod_{j=n+1,\ldots,n+m}^{\longleftarrow} \R_{ij}(z+u_i -v_{j-n}),\label{rnm12}\\
\cev{\R}_{nm}^{12}(u|v|z)= \prod_{i=1,\dots,n}^{\longleftarrow}
\prod_{j=n+1,\ldots,n+m}^{\longrightarrow} \R_{ij}(z+u_i -v_{j-n}).\label{rnm123}
\end{align}
In \eqref{rnm12} and \eqref{rnm123}, the superscripts $1$ and $2$ indicate the tensor factors in \eqref{nhk78} while the arrows indicate the order of the factors.
For example,
if $n=3, m=2$ and $\R_{ij}= \R_{ij}(z+u_i -v_{j-n})$
 we have
\beq\label{examplet}\R_{32}^{12}(u|v|z)=\R_{15}\R_{14}\R_{25}\R_{24}\R_{35}\R_{34} \fand\cev{\R}_{32}^{12}(u|v|z)=\R_{34}\R_{35}\R_{24}\R_{25}\R_{14} \R_{15}.\eeq
Observe that, due to the expansion convention  introduced in Section \ref{subsec011}, the expressions of the form $(z+u_i-v_{j-n})^r$ with $r<0$ are   expanded in  negative powers of the variable $z$, so that \eqref{rnm12} and \eqref{rnm123} contain only nonnegative powers of the variables $u_1,\ldots ,u_n$ and $v_{1},\ldots ,v_m$.
In order to simplify the notation, we write
\beq\label{rnm1234}
\R_{nm}^{12}(u|v)=\R_{nm}^{12}(u|v|0)\fand \cev{\R}_{nm}^{12}(u|v)=\cev{\R}_{nm}^{12}(u|v|0),
\eeq
where, due to the aforementioned expansion convention,  expressions of the form $( u_i-v_{j-n})^r$ with $r<0$  are    expanded in negative powers of the variable $u_i$, so that they contain only nonnegative powers of   $v_{1},\ldots ,v_m$.
The functions $R_{nm}^{12}(u|v|z)$, $\cev{R}_{nm}^{12}(u|v|z)$, $R_{nm}^{12}(u|v)$   and $\cev{R}_{nm}^{12}(u|v)$
corresponding to   $R$-matrix \eqref{R} can be defined analogously.

We will often combine the ordered product notation, as introduced in Section \ref{subsec021}, with the products of the form as in \eqref{rnm12} and \eqref{rnm123}. For example, in the expressions such as
$$\R_{nm}^{12}(u|v|z)\cdotlr X,\quad\text{ where}\quad X\in(\ndo\CC^N)^{\ot (n+m)},$$
 the  tensor factors $1,\ldots, n$ of $\R_{nm}^{12}(u|v|z)$ are applied from the left while the tensor factors $n+1,\ldots ,n+m$ are applied from the right, e.g., for $n=3$ and $m=2$ we have
$$\R_{32}^{12}(u|v|z)\cdotlr X =\R_{35}\cdotlr\left(\R_{34}\cdotlr\left(\R_{25}\cdotlr\left(\R_{24}\cdotlr\left(\R_{15}\cdotlr\left(\R_{14}\cdotlr X\right)\right)\right)\right)\right). $$

For any  integer $n\geqslant 1$, the variables $u=(u_1,\ldots, u_n)$ and the single variable $z$  define
\beq\label{teen}
T_{[n]}(u|z)=T_1(z+u_1)\ldots T_n(z+u_n )\fand T_{[n]}^+(u|z)=T^+_1(z+u_1)\ldots T^+_n(z+u_n).
\eeq
In particular, we write
$$
T_{[n]}(u)=T_1(u_1)\ldots T_n(u_n )\fand T_{[n]}^+(u)=T^+_1(u_1)\ldots T^+_n(u_n).
$$
We regard the coefficients of the series in \eqref{teen} as   operators on the vacuum  module $\Vc_c(\g_N)$. Hence the series  $T_{[n]}^+(u|z)$ belongs to $(\ndo\CC^N)^{\ot n}\ot \ndo\Vc_c(\g_N) [[u_1,\ldots ,u_n,z]]$.
Moreover, due to the expansion  convention introduced in Section \ref{subsec011} and   relations \eqref{DX3},  the series $T_{[n]}(u|z)$ belongs to
$(\ndo\CC^N)^{\ot n}\ot \om(\Vc_c(\g_N), \Vc_c(\g_N)[z^{-1}] [[u_1,\ldots ,u_n,h]])$.
As before, we use the arrows to indicate the opposite order of factors. For example,
$$
\cev{T}_{[n]}(u|z)= T_n(z+u_n)\ldots T_1(z+u_1) \fand \cev{T}_{[n]}^+(u|z)= T^+_n(z+u_n)\ldots T^+_1(z+u_1).
$$
As in \cite{EK}, by employing Yang--Baxter equation \eqref{YBE} and defining relations \eqref{DX1}--\eqref{DX3}   one obtains   the    more general form of the $RTT$ relations.

\begin{pro}\label{prop_rel}
 For any $c\in\CC$ and integers $n,m\geqslant 1$   the equalities
\begin{align}
&R_{nm}^{12}(u|v|z-w)T_{[n]}^{13}(u|z)T_{[m]}^{23}(v|w)
=T_{[m]}^{23}(v|w)T_{[n]}^{13}(u|z)R_{nm}^{12}(u|v|z-w),\label{rttz2}\\
&R_{nm}^{12}(u|v|z-w)T_{[n]}^{+13}(u|z)T_{[m]}^{+23}(v|w)
=T_{[m]}^{+23}(v|w)T_{[n]}^{+13}(u|z)R_{nm}^{12}(u|v|z-w),\label{rttz1}\\
&\R_{nm}^{ 12}(u|v|z-w+h \sigma c/2)T_{[n]}^{13}(u|z)T_{[m]}^{+23}(v|w)\non\\
&\qquad\qquad\qquad\qquad\qquad\qquad
=T_{[m]}^{+23}(v|w)T_{[n]}^{13}(u|z)\R_{nm}^{ 12}(u|v|z-w-h\sigma  c/2)
\label{rttz3}
\end{align}
hold for operators on
$
(\ndo\CC^{N})^{\ot  n} \ot
(\ndo\CC^{N})^{\ot  m}\ot  \Vc_c(\g_N).
$
\end{pro}

\section{$h$-adic quantum vertex algebras in types \texorpdfstring{$B,\,C\text{ and }D$}{B, C and D}
}\label{sec03}

\subsection{Vacuum module as an $h$-adic quantum vertex algebra}\label{subsec032}

The following notion of $h$-adic quantum vertex algebra was introduced by Li in \cite{Li}. As explained therein, it presents a slight generalization of the notion of quantum VOA, which was introduced by Etingof and Kazhdan in \cite{EK}.
From now on, the tensor products are understood as $h$-adically completed.
\begin{defn}\label{qvoa}
An {\em $h$-adic quantum vertex algebra} is a quadruple $(V,Y,\vac,\Sc)$ which satisfies the following axioms:
\begin{enumerate}[(1)]
\item  $V$ is a topologically free $\mathbb{C}[[h]]$-module.
\item $Y$ is the {\em vertex operator map}, a $\mathbb{C}[[h]]$-module map
\begin{align*}
Y \colon V\ot V&\to V((z))[[h]]\\
u\ot v&\mapsto Y(z)(u\ot v)=Y(u,z)v=\sum_{r\in\mathbb{Z}} u_r v \ts z^{-r-1}
\end{align*}
which satisfies the {\em weak associativity}:
for any $u,v,w\in V$ and $n\in\mathbb{Z}_{\geqslant 0}$
there exists $r\in\mathbb{Z}_{\geqslant 0}$
such that
\begin{equation}\label{associativity}
(z_0 +z_2)^r\ts Y(u,z_0 +z_2)Y(v,z_2)\ts w - (z_0 +z_2)^r\ts Y\big(Y(u,z_0)v,z_2\big)\ts w
\in h^n V[[z_0^{\pm 1},z_2^{\pm 1}]].
\end{equation}
\item $\vac$ is the {\em vacuum vector}, an element of $V$  which satisfies
\beq\label{v1}
Y(\vac ,z)v=v\quad\text{for all }v\in V,
\eeq
and for any $v\in V$ the series $Y(v,z)\ts\vac$ is a Taylor series in $z$ with
the property
\beq\label{v2}
\lim_{z\to 0} Y(v,z)\ts\vac =v.
\eeq
\item $\Sc=\Sc(z)$ is a $\mathbb{C}[[h]]$-module map
$ V\otimes V\to V\otimes V\otimes\mathbb{C}((z))$ which satisfies the {\em shift condition}
\begin{align}
&[\Dc\otimes 1, \mathcal{S}(z)]=-\frac{d}{dz}\mathcal{S}(z)\quad \text{for}\quad \Dc\in\ndo V\text{   defined by }\Dc v=v_{-2}\vac,\label{s1}\\
\intertext{the {\em Yang--Baxter equation}}
&\mathcal{S}_{12}(z_1)\ts\mathcal{S}_{13}(z_1+z_2)\ts\mathcal{S}_{23}(z_2)
=\mathcal{S}_{23}(z_2)\ts\mathcal{S}_{13}(z_1+z_2)\ts\mathcal{S}_{12}(z_1),\label{s2}\\
\intertext{the {\em unitarity condition}}
&\mathcal{S}_{21}(z)=\mathcal{S}^{-1}(-z),\label{s3}
\end{align}
the {\em $\mathcal{S}$-locality}:
for any $u,v\in V$ and $n\in\mathbb{Z}_{\geqslant 0}$ there exists
$r\in\mathbb{Z}_{\geqslant 0}$ such that
\begin{align}
&(z_1-z_2)^{r}\ts Y(z_1)\big(1\otimes Y(z_2)\big)\big(\mathcal{S}(z_1 -z_2)(u\otimes v)\otimes w\big)
\nonumber\\
&\quad-(z_1-z_2)^{r}\ts Y(z_2)\big(1\otimes Y(z_1)\big)(v\otimes u\otimes w)
\in h^n V[[z_1^{\pm 1},z_2^{\pm 1}]]\quad\text{for all }w\in V,\label{locality}
\end{align}
and the {\em hexagon identity}:
\begin{align}\label{hexagon}
\Sc(z_1)\left(Y(z_2)\ot 1\right) =\left(Y(z_2)\ot 1\right)\Sc_{23}(z_1)\Sc_{13}(z_1+z_2).
\end{align}
\end{enumerate}
\end{defn}

Denote by $\vac$ the image of the unit $1\in\DY_c(\g_N)$ in the vacuum module $\Vc_c(\g_N)$.
The next theorem  is a generalization of the Etingof--Kazhdan's construction of quantum VOAs in type $A$; see \cite[Theorem 2.3]{EK}.
\begin{thm}\label{EK:qva}
For any $c\in \CC$
there exists a unique  structure of $h$-adic quantum vertex algebra
on  $\V_c(\g_N)$, where $\g_N=\on_N,\spn_N$, such that the vacuum vector is
$\vac\in \V_c(\g_N)$, the vertex operator map is defined by
\beq\label{qva1}
Y\big(T_{[n]}^+ (u)\vac,z\big)=T_{[n]}^+ (u|z)\ts T_{[n]} (u|z+h\sigma c/2)^{-1}
\eeq
  and the map $\mathcal{S}(z)$ is defined by
\begin{align}
\mathcal{S}(z)\Big(\R_{nm}^{  12}(u|v|z)^{-1} \ts T_{[m]}^{+24}(v) \ts
\R	_{nm}^{  12}(u|v|z-h\sigma   c) \ts T_{[n]}^{+13}(u)\ts(\vac\otimes \vac) \Big)&\nonumber\\
=T_{[n]}^{+13}(u)  \ts\R_{nm}^{  12}(u|v|z+h\sigma   c)^{-1} \ts
 T_{[m]}^{+24}(v) \ts \R_{nm}^{  12}(u|v|z)\ts(\vac\otimes \vac)&\label{qva2}
\end{align}
for operators on
$
(\ndo\mathbb{C}^{N})^{\otimes n} \otimes
(\ndo\mathbb{C}^{N})^{\otimes m}\otimes \V_c(\g_N) \ot \V_c(\g_N)$.
\end{thm}

\begin{prf}
Recall that the $\CC[[h]]$-module $\V_c(\g_N)$ is topologically free.
Let us prove that the map $Y(z)$ is well-defined by \eqref{qva1}.
As the coefficients of matrix entries of all $T_{[n]}^+ (u)\vac$ span an $h$-adically dense subset of $\V_c(\g_N)$,
it is sufficient to show that $Y(z)$ maps the ideal of defining relations \eqref{DX2} and \eqref{TT}  for the dual Yangian $\Y^+(\g_N)$ to itself.
We only verify this for defining relations  \eqref{TT}.
As for    \eqref{DX2}, this follows by employing  Proposition \ref{prop_rel}, Yang--Baxter equation \eqref{YBE} and arguing as in the proof of \cite[Lemma 2.1]{EK}.

For any nonnegative integers $n,m$ and the variables $u=(u_1,\ldots ,u_n)$, $v=(v_1,\ldots ,v_m)$ and $w$
we apply   $Y(z)$ on the  expression
$$
a\coloneqq T_{[n]}^{+14}(u)\ts T_{[1]}^{+24}(w) \ts T_{[1]}^{+24}(w+h\kp)'\ts  T_{[m]}^{+34}(v)\vac.
$$
The coefficients of $a$ belong to the tensor product
$$
  \overbrace{(\ndo\CC^N)^{\ot n}}^{1}\ot \overbrace{\ndo\CC^N}^{2} \ot \overbrace{(\ndo\CC^N)^{\ot m}}^{3}\ot \overbrace{\V_c(\g_N)}^{4}
$$
and its superscripts indicate the tensor copies as indicated above.
Using \eqref{qva1}
we get
\begin{align}
&Y(a,z)=
T_{[n]}^{+14}(u|z)\ts
T_{[1]}^{+24}(w|z)\ts A\ts\ts
T_{[n]}^{14}(u|z+h\sigma c/2)^{-1},\qquad\text{where}\label{hfg1}\\
&\qquad A=
B^{24}\cdotlr\left(
 T_{[1]}^{+24}(w|z+h\kp)'\ts
T_{[m]}^{+34}(v|z)\ts
T_{[m]}^{34}(v|z+h\sigma c/2)^{-1}\right),\label{hfg2}\\
&\qquad B=T_{[1]}(z+w+h\sigma c/2)^{-1}\cdotlr \left(T_{[1]}(z+w+h\kp+h\sigma c/2)^{-1}\right)'.\non
\end{align}
Note that
\begin{align}
B'=T_{[1]}(z+w+h\kp+h\sigma c/2)^{-1}  \left(T_{[1]}(z+w+h\sigma c/2)^{-1} \right)'=1
\label{thesecondequality}
\end{align}
so, consequently, $B=1$. Indeed, the second equality in \eqref{thesecondequality} follows directly from
defining relations \eqref{DX1}.
Therefore, by combining \eqref{hfg1} and \eqref{hfg2} we find
\begin{align*}
Y(a,z)= &\,\,
T_{[n]}^{+14}(u|z)\ts
T_{[1]}^{+24}(w|z)\ts
 T_{[1]}^{+24}(w|z+h\kp)'\ts
T_{[m]}^{+34}(v|z)\\
&\cdot T_{[m]}^{34}(v|z+h\sigma c/2)^{-1}\ts
T_{[n]}^{14}(u|z+h\sigma c/2)^{-1}.
\end{align*}
Finally, due to \eqref{TT}  we have $T_{[1]}^{+24}(w|z)
 T_{[1]}^{+24}(w|z+h\kp)'=1$, so that
\begin{align*}
Y(a,z)=
T_{[n]}^{+14}(u|z)\ts
T_{[m]}^{+34}(v|z)\ts T_{[m]}^{34}(v|z+h\sigma c/2)^{-1}\ts
T_{[n]}^{14}(u|z+h\sigma c/2)^{-1}.
\end{align*}
It is  now clear that the expression $Y(a,z)$ coincides with the image of $T_{[n]}^{+14}(u)  T_{[m]}^{+34}(v)\vac$ with respect to \eqref{qva1}, so we conclude that the map $Y(z)$ is well-defined.

Next, we  prove that the map $\Sc(z)$ is well-defined.
Due to crossing symmetry property \eqref{csym2}, we express \eqref{qva2} as
\begin{align}
\mathcal{S}(z)\big(
   T_{[n]}^{+13}(u)\ts T_{[m]}^{+24}(v) \ts (\vac\otimes \vac) \big)
=&\,\,\cev\Rt_{nm}^{  12}(u|v|z-h\kp-h\sigma   c)\cdotlr \Big(\R_{nm}^{  12}(u|v|z) T_{[n]}^{+13}(u)\Big. \label{mapes}\\
\Big.   &\cdot\R_{nm}^{  12}(u|v|z+h\sigma   c)^{-1} \ts
 T_{[m]}^{+24}(v) \ts \R_{nm}^{  12}(u|v|z)\ts(\vac\otimes \vac) \hspace{-2pt}\Big),\non
\end{align}
where, in consistency with notation introduced in \eqref{rnm12} and \eqref{rnm123}, we write
\begin{align*}
\Rt_{nm}^{12}(u|v|z)= \prod_{i=1,\dots,n}^{\longrightarrow}
\prod_{j=n+1,\ldots,n+m}^{\longleftarrow} \R_{ij}(z+u_i -v_{j-n})',\\
\cev{\Rt}_{nm}^{12}(u|v|z)= \prod_{i=1,\dots,n}^{\longleftarrow}
\prod_{j=n+1,\ldots,n+m}^{\longrightarrow} \R_{ij}(z+u_i -v_{j-n})'.
\end{align*}
Therefore, the map $\Sc(z)$ coincides with the composition $S^{(1)} (z)\circ S^{(2)}(z)\circ S^{(3)} (z)\circ S^{(4)} (z)$ of the following maps:
\begin{align*}
&
 T_{[n]}^{+13}(u)\ts T_{[m]}^{+24}(v)\ts (\vac\otimes \vac)
\,\xmapsto{S^{(1)}(z)}\,
T_{[n]}^{+13}(u)\ts\R_{nm}^{  12}(u|v|z+h\sigma c)^{-1}
\ts T_{[m]}^{+24}(v)
\ts (\vac\otimes \vac),\\
& T_{[n]}^{+13}(u)\ts T_{[m]}^{+24}(v)\ts (\vac\otimes \vac)
\,\xmapsto{S^{(2)}(z)}\,
T_{[n]}^{+13}(u)\ts T_{[m]}^{+24}(v)
\ts\R_{nm}^{  12}(u|v|z)
\ts (\vac\otimes \vac), \\
&
 T_{[n]}^{+13}(u)\ts T_{[m]}^{+24}(v)\ts (\vac\otimes \vac)
\,\xmapsto{S^{(3)}(z)}\,
\R_{nm}^{  12}(u|v|z)
\ts
T_{[n]}^{+13}(u)\ts T_{[m]}^{+24}(v)
\ts(\vac\otimes \vac), \\
&
 T_{[n]}^{+13}(u)\ts T_{[m]}^{+24}(v)\ts (\vac\otimes \vac)
\,\xmapsto{S^{(4)}(z)}\,
T_{[m]}^{+24}(v) \ts
\cev{\Rt}_{nm}^{  12}(u|v|z-h\kp-h\sigma   c)
\ts
T_{[n]}^{+13}(u)
\ts(\vac\otimes \vac).
\end{align*}
Hence it is sufficient to check that all  $S^{(i)}(z)$ are well-defined maps on $\Vc_c(\g_N)$, i.e. that they map the ideal of defining relations \eqref{DX2} and \eqref{TT} for the dual Yangian $\Y^+(\g_N)$ to itself. The fact that the given maps preserve relation \eqref{DX2} can be proved  by using
  Yang--Baxter equation \eqref{YBE} and arguing as in the proof of \cite[Lemma 2.1]{EK}. As for relation \eqref{TT}, this follows  by employing
  crossing symmetry properties \eqref{csym} and \eqref{csym2}.
	
The   weak associativity \eqref{associativity}, Yang--Baxter equation \eqref{s2}, unitarity property \eqref{s3} and $\Sc$-locality property \eqref{locality} can be verified by straightforward calculations    which rely on the properties of the $R$-matrix $\R(u)$ provided  in Section \ref{subsec021} and on Proposition \ref{prop_rel}. They closely follow the corresponding proofs in type $A$, as given in \cite[Theorem 4.1]{JKMY}.	
Regarding the axioms concerning the vacuum vector $\vac$, \eqref{v1} is clear and \eqref{v2} is a consequence of the identity $T(u)\vac =\vac$.

Let us prove shift condition \eqref{s1}. Let $n,m\geqslant 0$ be arbitrary integers. By applying \eqref{qva1} on $\vac$ and then taking the coefficient of the variable $z$ we obtain
\beq\label{Doperator}
\Dc  \ts T_{[n]}^{+}(u)\vac = \textstyle\left(\sum_{k=1}^n  \frac{\partial}{\partial u_k}\right)T_{[n]}^{+}(u)\vac.
\eeq
In particular, note that \eqref{v1} implies $\Dc\vac=0$.
Therefore, by applying $\Dc \ot 1$ on \eqref{mapes} we find that
$
\left(\Dc \ot 1\right)\mathcal{S}(z)\big(
   T_{[n]}^{+13}(u)\ts T_{[m]}^{+24}(v) \ts (\vac\otimes \vac) \big)
	$
	is equal to
	\begin{align}
\cev\Rt_{nm}^{  12}(u|v|z-h\kp-h\sigma   c)\cdotlr \Big(\R_{nm}^{  12}(u|v|z) \textstyle\left( \left(\sum_{k=1}^n  \frac{\partial}{\partial u_k}\right)T_{[n]}^{+13}(u)\right)\Big. \non\\
\Big.   \cdot \R_{nm}^{  12}(u|v|z+h\sigma   c)^{-1} \ts
 T_{[m]}^{+24}(v) \ts \R_{nm}^{  12}(u|v|z)\ts(\vac\otimes \vac) \Big).\label{s1proof1}
\end{align}
On the other hand, by applying the map $\Sc(z)$ on
$$\left(\Dc \ot 1\right)\big(
   T_{[n]}^{+13}(u)\ts T_{[m]}^{+24}(v) \ts (\vac\otimes \vac) \big)=\textstyle\left( \left(\sum_{k=1}^n  \frac{\partial}{\partial u_k}\right)T_{[n]}^{+13}(u)\right)\ts T_{[m]}^{+24}(v) \ts (\vac\otimes \vac) $$
	we get
		\begin{align}
		\textstyle  \left(\sum_{k=1}^n  \frac{\partial}{\partial u_k}\right)\bigg(
\cev\Rt_{nm}^{  12}(u|v|z-h\kp-h\sigma   c)\cdotlr \Big(\R_{nm}^{  12}(u|v|z)\ts T_{[n]}^{+13}(u)\Big.\bigg. \non\\
\bigg. \Big.   \cdot \R_{nm}^{  12}(u|v|z+h\sigma   c)^{-1} \ts
 T_{[m]}^{+24}(v) \ts \R_{nm}^{  12}(u|v|z)\ts(\vac\otimes \vac) \Big)\bigg).\label{s1proof2}
\end{align}
Finally,   using the observation $\textstyle \frac{\partial}{\partial u_k}\R(u_k+z) =  \frac{\partial}{\partial z}\R(u_k+z)$ we see that  the difference of expressions \eqref{s1proof1} and \eqref{s1proof2} coincides with
$$- \frac{\partial}{\partial z} \Sc(z)\big(
   T_{[n]}^{+13}(u)\ts T_{[m]}^{+24}(v) \ts (\vac\otimes \vac) \big), $$
so that shift condition \eqref{s1} follows.

It remains to verify hexagon identity \eqref{hexagon}.
Its proof goes similarly as the proof in type $A$; see proof of \cite[Theorem 2.3.8]{G}.
 Let $n,m,k$ be arbitrary nonnegative integers. Label the tensor copies as follows
\beq\label{hex1}
\overbrace{\left(\ndo\CC^N\right)^{\ot n}}^{1} \ot
\overbrace{\left(\ndo\CC^N\right)^{\ot m}}^{2} \ot
\overbrace{\left(\ndo\CC^N\right)^{\ot k}}^{3} \ot
\overbrace{\Vc_c(\g_N)}^{4}\ot\overbrace{\Vc_c(\g_N)}^{5}\ot\overbrace{\Vc_c(\g_N)}^{6}.
\eeq
First,  we apply the left hand side of hexagon identity \eqref{hexagon} on the expression
\beq\label{hex22}
T_{[n]}^{+14}(u)\ts T_{[m]}^{+25}(v)\ts T_{[k]}^{+36}(w)\ts(\vac\ot\vac\ot\vac),
\eeq
whose coefficients, with respect to the variables $u=(u_1,\ldots ,u_n)$, $v=(v_1,\ldots ,v_m)$ and $w=(w_1,\ldots ,w_k)$,  belong to \eqref{hex1}. By applying $Y(z_2)\ot 1$ on \eqref{hex22} we obtain the expression
\beq\label{hex2}
T_{[n]}^{+14}(u|z_2)\ts T_{[n]}^{14}(u|z_2 +h\sigma c/2)^{-1}\ts T_{[m]}^{+24}(v)\ts T_{[k]}^{+35}(w)\ts(\vac\ot\vac).
\eeq
By combining relation \eqref{rttz3} with crossing symmetry properties \eqref{csym2} and \eqref{csym3}  we rewrite \eqref{hex2} as
\begin{gather}
A\cdotrl \left( T_{[n]}^{+14}(u|z_2)\ts T_{[m]}^{+24}(v)\ts T_{[k]}^{+35}(w)\ts B\ts \right) (\vac\ot\vac),\qquad\text{where}\label{hex3}\\
A=\R_{nm}^{12}(u|v|z_2+h\sigma c+2h\kp)  \Fand B=\R_{nm}^{12}(u|v|z_2)^{-1} .\non
\end{gather}
Due to   \eqref{mapes}, by applying the map $\Sc(z_1)$ on \eqref{hex3} we get
\beq\label{hex4}
A\cdotrl \left(H\cdotlr \left(K\ts T_{[n]}^{+14}(u|z_2)\ts  T_{[m]}^{+24}(v)\ts L\ts  T_{[k]}^{+35}(w)\ts K\ts B\right)\right)  (\vac\ot\vac),
\eeq
where
\begin{alignat*}{3}
&H=H_{[2]}H_{[1]},\quad &&H_{[1]}= \cev{\Rt}_{nk}^{13}(u|w|z_1+z_2-h\sigma c-h\kp),\quad &&H_{[2]}=\cev{\Rt}_{mk}^{23}(v|w|z_1-h\sigma c-h\kp)  ,\\
&K=K_{[1]}K_{[2]},\quad  && K_{[1]}=\R_{nk}^{13}(u|w|z_1 +z_2),\quad &&K_{[2]}= \R_{mk}^{23}(v|w|z_1),\\
&L=L_{[2]}L_{[1]} ,\quad &&L_{[1]}=\R_{nk}^{13}(u|w|z_1 +z_2 +h\sigma c)^{-1},\quad &&L_{[2]}=\R_{mk}^{23}(v|w|z_1 +h\sigma c)^{-1}.
\end{alignat*}

We  now apply the right hand side of hexagon identity \eqref{hexagon} on \eqref{hex22} and show that the result coincides with \eqref{hex4}.
The tensor copies are again labelled as in \eqref{hex1}.
First, applying the map $\Sc_{46}(z_1 +z_2)$ on \eqref{hex22} we obtain
\beq\label{hexx}
H_{[1]}\cdotlr\left(K_{[1]}\ts T_{[n]}^{+14}(u)\ts L_{[1]}\ts  T_{[m]}^{+25}(v)\ts T_{[k]}^{+36}(w)\ts K_{[1]}\right)(\vac\ot\vac\ot\vac).
\eeq
Next,   applying the map $\Sc_{56}(z_1)$ on \eqref{hexx} we get
\begin{align*}
&H_{[1]}\cdotlr\left(K_{[1]}\ts T_{[n]}^{+14}(u)\ts L_{[1]}\left(
H_{[2]}\cdotlr\left( K_{[2]}\ts
T_{[m]}^{+25}(v)\ts L_{[2]}\ts T_{[k]}^{+36}(w)\ts K_{[2]}\right)\right) K_{[1]}\right)(\vac\ot\vac\ot\vac)\\
=\,\,&H_{[1]}\cdotlr\left(
H_{[2]}\cdotlr\left(  K_{[1]}\ts T_{[n]}^{+14}(u)\ts L_{[1]} \ts K_{[2]}\ts
T_{[m]}^{+25}(v)\ts L_{[2]}\ts T_{[k]}^{+36}(w)\ts K_{[2]}\right) K_{[1]}\right)(\vac\ot\vac\ot\vac).
\end{align*}
Finally, by applying $Y(z_2)\ot 1$ we get
\begin{align}
H_{[1]}\cdotlr\left(
H_{[2]}\cdotlr\left(  K_{[1]}\ts  X\ts L_{[2]}\ts T_{[k]}^{+35}(w)\ts K_{[2]}\right)  K_{[1]}\right) (\vac\ot\vac ),\label{hex9}
\end{align}
where
\begin{align*}
X=& \,\,T_{[n]}^{+14}(u|z_2)\ts T_{[n]}^{14}(u|z_2+h\sigma c/2)^{-1} \ts L_{[1]} \ts K_{[2]}\ts
T_{[m]}^{+24}(v)\\
=&\,\, L_{[1]}\cdotrl\left(  K_{[2]}\ts  T_{[n]}^{+14}(u|z_2)\ts T_{[n]}^{14}(u|z_2+h\sigma c/2)^{-1}\ts
T_{[m]}^{+24}(v)\right).
\end{align*}
As before, we combine relation \eqref{rttz3} with crossing symmetry properties \eqref{csym2} and \eqref{csym3} to write $X$ as
\beq\label{hex8}
X= \ts L_{[1]}\cdotrl\left(  K_{[2]}\ts  T_{[n]}^{+14}(u|z_2)\left( A\cdotrl \left(
T_{[m]}^{+24}(v)   \ts B\ts T_{[n]}^{14}(u|z_2+h\sigma c/2)^{-1}\right)\right)\right).
\eeq
Identities  $T (z)\vac =\vac$ and \eqref{hex8}  imply that   \eqref{hex9}  is equal to
\begin{align}
H_{[1]}\cdotlr\left(
H_{[2]}\cdotlr\left(  K_{[1]}\ts  Z\ts L_{[2]}\ts T_{[k]}^{+35}(w)\ts K_{[2]}\right)  K_{[1]}\right) (\vac\ot\vac ),\label{hexa}
\end{align}
where
\beq\label{hexb}
Z= \ts L_{[1]}\cdotrl\left(  K_{[2]}\ts  T_{[n]}^{+14}(u|z_2)\left( A\cdotrl \left(
T_{[m]}^{+24}(v)   \ts B\ts \right)\right) \right).
\eeq

The following consequences of Yang--Baxter equation \eqref{YBE} and crossing symmetry property \eqref{csym3} can be   verified by a straightforward calculation:
\begin{alignat}{2}
&L_{[1]}\cdotrl\left(  K_{[2]}\ts   A\right)   = A\cdotrl \left(K_{[2]}\ts L_{[1]} \right), \quad
&&H_{[2]} \ts A \ts K_{[1]}= K_{[1]}\ts A\ts H_{[2]},\label{hexxq2}\\
&B  \ts L_{[1]}      \ts L_{[2]}=L_{[2]}  \ts L_{[1]}      \ts B=L B,\qquad    &&B\ts K_{[2]} \ts K_{[1]} = K_{[1]}\ts K_{[2]} \ts B= K B.\label{hexxq3}
\end{alignat}
By using the first equality in \eqref{hexxq2}  we can  write \eqref{hexb} as
$$
Z=  T_{[n]}^{+14}(u|z_2) \left( A\cdotrl \left(K_{[2]} \ts T_{[m]}^{+24}(v)\ts B  \ts L_{[1]}      \right)\right).
$$
Therefore, the original expression in \eqref{hexa} is equal to
\begin{align*}
H_{[1]}\cdotlr\left(
H_{[2]}\cdotlr\left(  K_{[1]}\ts
T_{[n]}^{+14}(u|z_2) \left( A\cdotrl \left(K_{[2]} \ts T_{[m]}^{+24}(v)\ts B  \ts L_{[1]}      \right)\right) L_{[2]}\ts T_{[k]}^{+35}(w)\ts K_{[2]}\right)  K_{[1]}\right)(\vac\ot\vac ).
\end{align*}
By moving the element $A$ to the left and using $K=K_{[1]}K_{[2]}$ we obtain
\begin{align*}
H_{[1]}\cdotlr\left(
H_{[2]}\cdotlr\left(
A\cdotrl \left(
 K_{[1]}\ts
T_{[n]}^{+14}(u|z_2)\ts   K_{[2]} \ts T_{[m]}^{+24}(v)\ts B  \ts L_{[1]}      \ts L_{[2]}\ts T_{[k]}^{+35}(w)\ts K_{[2]}\right)\right)  K_{[1]}\right)(\vac\ot\vac )&\\
=H_{[1]}\cdotlr\left(
H_{[2]}\cdotlr\left(
A\cdotrl \left(
 K \ts
T_{[n]}^{+14}(u|z_2)  \ts T_{[m]}^{+24}(v)\ts B  \ts L_{[1]}      \ts L_{[2]} \ts T_{[k]}^{+35}(w)\ts K_{[2]}\right)\right)  K_{[1]}\right)(\vac\ot\vac )&.
\end{align*}
Next, we employ the second equality in \eqref{hexxq2} and $H=H_{[2]}H_{[1]}$ to write the given expression as
\begin{align}
&A\cdotrl \left(
H_{[1]}\cdotlr\left(
H_{[2]}\cdotlr\left(
 K \ts
T_{[n]}^{+14}(u|z_2)  \ts T_{[m]}^{+24}(v)\ts B  \ts L_{[1]}      \ts L_{[2]} \ts T_{[k]}^{+35}(w)\ts K_{[2]} \ts K_{[1]}\right)\right)\right)(\vac\ot\vac )\non\\
=\,\,&A\cdotrl \left(
H\cdotlr\left(
 K \ts
T_{[n]}^{+14}(u|z_2)  \ts T_{[m]}^{+24}(v)\ts B  \ts L_{[1]}      \ts L_{[2]} \ts T_{[k]}^{+35}(w)\ts K_{[2]} \ts K_{[1]}\right)\right)(\vac\ot\vac ). \label{hexg}
\end{align}
Finally, we use both equalities in \eqref{hexxq3}  to move the element $B$ in \eqref{hexg} to the right, thus getting \eqref{hex4}, as required. Therefore, we conclude that hexagon identity \eqref{hexagon} holds, so  the proof of the theorem is over.
\end{prf}

\subsection{Center of the $h$-adic quantum vertex algebra at the critical level }\label{subsec042}
In this section, we consider the $h$-adic quantum vertex algebra $\Vcr(\g_N)= \Vc_{c_{\textrm{crit}}}(\g_N)$ at the critical level
$$
c_{\textrm{crit}} =-\frac{2\kp}{\sigma}=\begin{cases}
-N+2,&\text{if }\g_N=\on_N,\\
-\frac{N}{2}-1,&\text{if }\g_N=\spn_N.
\end{cases}
$$
First, we follow the exposition in \cite{M} to recall a particular case of the fusion procedure for the Brauer algebra  \cite{B}; see also \cite[Section 1.2]{Mnew}.
Let $\omega$ be an indeterminate and let $\Bc_m(\omega)$ be the Brauer algebra over the field $\CC(\omega)$ generated by the elements $s_1,\ldots ,s_{m-1}$ and $\varepsilon_1,\ldots ,\varepsilon_{m-1}$ subject to the defining relations as in \cite[Section 3]{M}.
Its complex subalgebra  generated by the elements $s_1,\ldots ,s_{m-1}$ is isomorphic to the group algebra of the symmetric group $\mathfrak{S}_m$, so that the elements $s_i$ are identified with the transpositions $(i,i+1)$. Let $s_{ij}\in \Bc_m(\omega)$, $i<j$, be the element corresponding to the transposition $(i,j)$ with respect to that isomorphism. Introduce  $\varepsilon_{ij}\in \Bc_m(\omega)$ by
$\varepsilon_{j-1\ts j}=\varepsilon_{j-1}$ and $\varepsilon_{ij}= s_{i\ts j-1}\varepsilon_{j-1}s_{i\ts j-1}$ for $i<j-1$.
Let $s^{(m)}\in\Bc_m(\omega)$ be the idempotent corresponding to the one-dimensional representation of the Brauer algebra which maps all $s_{ij}$ to the identity operator and all $\varepsilon_{ij} $ to the zero operator. It satisfies
\beq\label{important}
s_{ij}s^{(m)}=s^{(m)} s_{ij}=s^{(m)}\fand \varepsilon_{ij}s^{(m)}=s^{(m)} \varepsilon_{ij}=0.
\eeq

Consider the expression
$$
R(u_1,\ldots ,u_m) =\frac{1}{m!}\prod_{1\leqslant i<j\leqslant m}R_{ij}(u_i -u_j),
$$
where the products are taken in the lexicographical order on the pairs $(i,j)$. Define
 \beq\label{uem1}
u_{[m]}=\begin{cases}
(u-(m-1)h,\ldots,u-h, u)&\text{for }\g_N=\on_N\text{ and }m=1,\ldots ,N,\\
(u,u-h,\ldots, u-(m-1)h)&\text{for }\g_N=\spn_N\text{ and }m=1,\ldots ,N/2.
\end{cases}
\eeq
In the orthogonal case, let $S_{[m]}$ with $m=1,\ldots ,N$ denote the action of the idempotent $s^{(m)}\in \Bc_m (N)$ on the tensor product space
$(\CC^N)^{\ot m}$ with respect to the representation defined by
$s_{ij}\mapsto P_{ij}$ and $\varepsilon_{ij}\mapsto Q_{ij}$ for $i<j$. In the symplectic case, let
$S_{[m]}$ with $m=1,\ldots ,N/2$ denote the action of the idempotent $s^{(m)}\in\Bc_m (-N)$ on the tensor product space
$(\CC^N)^{\ot m}$ with respect to the representation defined by
$s_{ij}\mapsto -P_{ij}$ and $\varepsilon_{ij}\mapsto -Q_{ij}$ for $i<j$.
Due to a particular case of the  fusion procedure   for the Brauer algebra $\Bc_m (\omega)$, see \cite{IM,IMO},  we have
\beq\label{fusion}
 S_{[m]}=R(u_{[m]}) .
\eeq

Consider the tensor product
\beq\label{exprfsion3}
\overbrace{\ndo\CC^N}^{0} \ot \overbrace{(\ndo\CC^N)^{\ot m}}^{1} \ot\overbrace{\Vcr(\g_N)}^{2}.
\eeq
We use the following consequences of the fusion procedure    in the proof of Theorem \ref{centerthm}.

\begin{lem}\label{lemma9}
For any  $m=1,\ldots , N$ in the orthogonal case, $m=1,\ldots, N/2$ in the symplectic case and $\alpha\in \CC$   we have
\begin{align}
S_{[m]}^1 \ts \R_{1m}^{01}(v+h\alpha |u_{[m]})&=\cev{\R}_{1m}^{01}(v+h\alpha |u_{[m]})\ts S_{[m]}^1,\label{fusion3}\\
S_{[m]}^1 \ts T_{[m]}^{12} (u_{[m]}+h\alpha) &= \cev{T}_{[m]}^{12} (u_{[m]}+h\alpha) \ts S_{[m]}^1,\label{fusion22}
\\
S_{[m]}^1 \ts T_{[m]}^{+12} (u_{[m]}) &= \cev{T}_{[m]}^{+12} (u_{[m]}) \ts S_{[m]}^1\label{fusion2},
\end{align}
where the superscripts indicate the tensor factors in
\eqref{exprfsion3}.
\end{lem}

\begin{prf}
Yang--Baxter equation \eqref{YBE} implies
\beq\label{fusionprf2}
R(u_1,\ldots ,u_m)\ts \R_{1m}^{01}(v+h\alpha |u)=\cev{\R}_{1m}^{01}(v+h\alpha |u)\ts R(u_1,\ldots ,u_m)
\eeq
and defining relations \eqref{DX1} and \eqref{DX2} imply
\begin{align}
R(u_1,\ldots ,u_m)\ts T_{[m]}^{12} (u+h\alpha )&=\cev{T}_{[m]}^{12} (u+h\alpha )\ts R(u_1,\ldots ,u_m)\\
R(u_1,\ldots ,u_m)\ts T_{[m]}^{+12} (u )&=\cev{T}_{[m]}^{+12} (u )\ts R(u_1,\ldots ,u_m)\label{fusionprf}
\end{align}
for the variables $u=(u_1,\ldots , u_m)$,
where   $R(u_1,\ldots ,u_m)$ is applied on the tensor factor $(\ndo\CC^N)^{\ot m}$   of \eqref{exprfsion3} and $u+h\alpha=(u_1+h\alpha,\ldots , u_m+h\alpha)$.
 By evaluating the variables $u=(u_1,\ldots , u_m)$ at $u_{[m]}$ in equalities \eqref{fusionprf2}--\eqref{fusionprf}   and then applying fusion procedure \eqref{fusion} we obtain \eqref{fusion3}--\eqref{fusion2}, as required.
\end{prf}

Let $V$ be an $h$-adic quantum vertex algebra. As in \cite{JKMY}, we define the {\em center} of $V$   in analogy with vertex algebra theory, see, e.g., \cite[Chapter 3.3]{F}, as a $\CC[[h]]$-submodule
$$
\z(V)=\left\{v\in V\, :\, Y(w,z)v \in V[[z]] \text{ for all }w\in V\right\}.
$$
Due to \eqref{qva1}, the center of $\Vc_c(\g_N)$ coincides with the  $\Y(\g_N)$-invariants, i.e.
\beq\label{centar}
\z(\Vc_c(\g_N))=\left\{v\in\Vc_c(\g_N)\,:\, t_{ij}(u)v=\delta_{ij}v\text{ for all }i,j=1,\ldots ,N \right\}.
\eeq
As  in \cite[Chapter 11.2]{Mnew}, define the series
\beq\label{tcenter}
\TT_{m}^+ (u) =\tr_{1,\ldots,m} \ts S_{[m]} \ts T_{[m]}^+ (u_{[m]})\vac
= \tr_{1,\ldots,m} \ts S_{[m]} \ts T_{1}^+ (u_{1})\ldots  T_{m}^+ (u_{m})\vac
\eeq
in $\Vcr(\g_N)[[u]]$,
where $m=1,\ldots , N$ in the orthogonal case,  $m=1,\ldots, N/2$ in the symplectic case and the trace is taken over all $m$ copies of $\ndo\CC^N$.

\begin{thm}\label{centerthm}
All coefficients of $\TT_{m}^+ (u)$ belong to the center of the $h$-adic quantum vertex algebra $\Vcr(\g_N)$.
\end{thm}

\begin{prf}
Recall \eqref{uem1}. Consider the expressions
\beq\label{exprfsion}
A=\R_{1m}^{01}(v-h\kp|u_{[m]})^{-1}  \fand  B=\R_{1m}^{01}(v+h\kp|u_{[m]}).
\eeq
Their coefficients with respect to the variables $v$ and $u$ belong to  tensor product
\eqref{exprfsion3}
and their superscripts indicate the tensor copies therein.
By \eqref{fusion3} we have
\beq\label{exprfsion2}
S_{[m]}^1  A =\cev{A} \ts S_{[m]}^1\fand S_{[m]}^1\ts B =\cev{B} \ts S_{[m]}^1 ,
\eeq
where the superscript $1$ indicates that the idempotent $S_{[m]}$ acts on the tensor copy $(\ndo\CC^N)^{\ot m}$ of \eqref{exprfsion3}. Crossing symmetry property \eqref{csym} implies
$$A=\left(\R_{1m}^{01}(v|u_{[m]})\right)'^{{}_0}=\cev{\Rt}_{1m}^{01}(v|u_{[m]}) ,$$
where the transposition $'$ in the middle term is applied  on the tensor copy $\ndo\CC^N$ of \eqref{exprfsion3}.
Therefore, by crossing symmetry property \eqref{csym2} we have
\beq\label{exprfsion4}
A \cdotlr B =1.
\eeq

In order to prove the  theorem, it is sufficient to verify that
the coefficients of $\TT_{m}^+ (u)$ belong to the $\CC[[h]]$-submodule of $\Y(\g_N)$-invariants \eqref{centar}, i.e.
\beq\label{exprfsion5}
T(v) \ts\TT_{m}^+ (u)=\TT_{m}^+ (u).
\eeq
By using \eqref{tcenter} we write the left hand side in \eqref{exprfsion5}  as
\beq\label{fsn2}
\tr_{1,\ldots,m}\ts T^{02}_{[1]} (v)\ts S_{[m]}^1 \ts T_{[m]}^{+12} (u_{[m]}),
\eeq
where the superscripts indicate the tensor factors  in \eqref{exprfsion3}. Next, by \eqref{DX3} and $T(v)\vac =\vac$ we conclude that \eqref{fsn2} is equal to
\beq\label{fsn3}
\tr_{1,\ldots,m} \ts S_{[m]}^1  A\ts T_{[m]}^{+12} (u_{[m]})\ts B.
\eeq
We now  combine the cyclic property of the trace and $(S_{[m]})^2 =S_{[m]}$ with equalities \eqref{fusion2} and \eqref{exprfsion2}    to rewrite \eqref{fsn3} as follows:
\begin{align}
&\tr_{1,\ldots,m} \ts S_{[m]}^1 \ts A \ts T_{[m]}^{+12} (u_{[m]})\ts B
=
\tr_{1,\ldots,m} \ts  \cev{A} \ts S_{[m]}^1\ts T_{[m]}^{+12} (u_{[m]})\ts B \non\\
 =\,&\tr_{1,\ldots,m} \ts  \cev{A} \ts (S_{[m]}^1)^2 \ts T_{[m]}^{+12} (u_{[m]})\ts B
=\tr_{1,\ldots,m} \ts S_{[m]}^1 \ts A  \ts \cev{T}_{[m]}^{+12} (u_{[m]})\ts \cev{B} \ts S_{[m]}^1\non\\
  =\,&\tr_{1,\ldots,m}   \ts A \ts \cev{T}_{[m]}^{+12} (u_{[m]})\ts \cev{B} \ts (S_{[m]}^1)^2
=\tr_{1,\ldots,m}   \ts A  \ts \cev{T}_{[m]}^{+12} (u_{[m]})\ts \cev{B} \ts S_{[m]}^1\non\\
  =\,&\tr_{1,\ldots,m}   \ts A \ts S_{[m]}^1 \ts T_{[m]}^{+12} (u_{[m]})\ts B .\label{fsn6}
\end{align}
Finally, by the cyclic property of the trace,  expression \eqref{fsn6} equals
\begin{align*}
\tr_{1,\ldots,m}   \ts A \cdotlr\left( S_{[m]}^1 \ts T_{[m]}^{+12} (u_{[m]})\ts B \right)
=\tr_{1,\ldots,m}   \ts  S_{[m]}^1 \ts T_{[m]}^{+12} (u_{[m]})\left( A \cdotlr B  \right).
\end{align*}
By \eqref{exprfsion4}  this is equal to  $\TT_{m}^+ (u)$, so equality \eqref{exprfsion5} follows and the proof is over.
\end{prf}

In the orthogonal case, the series $\TT_{m}^+ (u)$ can be also written in the form
\beq\label{intheform}
\TT_{m}^+ (u) = \tr_{1,\ldots ,m}\, S_{[m]}\ts T_1^+ (u)\ts T_2^+(u-h)\ldots T_m^+(u-(m-1)h).
\eeq
Indeed, by the first equality in \eqref{important} we have $P_{ij}S_{[m]}P_{ij}=S_{[m]}$ for any $i,j=1,\ldots ,m$.
Hence conjugating   the expression under the trace by a suitable element  of the symmetric group $\mathfrak{S}_m$ we rewrite \eqref{intheform} as
$$
\tr_{1,\ldots ,m}\, S_{[m]}\ts T_m^+ (u)\ts T_{m-1}^+(u-h)\ldots T_1^+(u-(m-1)h).
$$
Next, using the cyclic property of the trace we move the idempotent $S_{[m]}$ to the right thus getting
 $$
\tr_{1,\ldots ,m}\,   T_m^+ (u)\ts T_{m-1}^+(u-h)\ldots T_1^+(u-(m-1)h) \ts S_{[m]}.
$$
Finally, fusion procedure  \eqref{fusion2} implies that this equals   $\TT_{m}^+ (u)$, as required.

\begin{thm}\label{thmmain4}
In $\Vcr(\g_N)$ we have
\beq\label{fixedpt}
\Sc(z) (\TT_{k}^+ (u)\ot \TT_{m}^+ (v)) =\TT_{k}^+ (u)\ot \TT_{m}^+ (v).
\eeq
\end{thm}

\begin{prf}
Label the tensor copies as follows:
\beq\label{exprfsion378}
\overbrace{(\ndo\CC^N)^{\ot k}}^{1} \ot \overbrace{(\ndo\CC^N)^{\ot m}}^{2} \ot\overbrace{\Vcr(\g_N)}^{3}\ot\overbrace{\Vcr(\g_N)}^{4}.
\eeq
Let
$$
A=\cev\Rt_{km}^{  12}(u_{[k]}|v_{[m]}|z+h\kp), \quad
B= \R_{km}^{  12}(u_{[k]}|v_{[m]}|z), \quad
C=\R_{km}^{  12}(u_{[k]}|v_{[m]}|z-2h\kp)^{-1},
$$
where
$
u_{[k]}=(\wvr{u}_1,\ldots ,\wvr{u}_{k})$ and $v_{[m]}=(\wvr{v}_1,\ldots ,\wvr{v}_m)$
are defined by \eqref{uem1}.
Using  crossing symmetry properties \eqref{csym}, \eqref{csym2} and \eqref{csym3} we find
\beq\label{csyms}
B\ts A=1\fand B\cdotrl C=1.
\eeq
Set
\begin{align*}
\wht{A}&=\prod_{i=1,\ldots ,k}^{\longrightarrow}\prod_{j=k+1,\ldots ,k+m}^{\longrightarrow}
\R_{ij}(z+\wvr{u}_i-\wvr{v}_{j-k}+h\kp)',\\
\wht{B}&=\prod_{i=1,\ldots ,k}^{\longleftarrow}\prod_{j=k+1,\ldots ,k+m}^{\longleftarrow}
\R_{ij}(z+\wvr{u}_i-\wvr{v}_{j-k}),\\
\wht{C}&=\prod_{i=1,\ldots ,k}^{\longrightarrow}\prod_{j=k+1,\ldots ,k+m}^{\longrightarrow}
\R_{ij}(z+\wvr{u}_i-\wvr{v}_{j-k}-2h\kp)^{-1}.
\end{align*}
By using fusion procedure \eqref{fusion} and arguing as in the proof of Lemma \ref{lemma9} one can prove
\beq\label{zabc}
S_{[k]}^1 A=\wht{A}S_{[k]}^1,\qquad  S_{[k]}^1 B=\wht{B}S_{[k]}^1,\Fand S_{[k]}^1 C=\wht{C}S_{[k]}^1.
\eeq

Explicit formula \eqref{mapes} for   the operator $\Sc(z)$ at the critical level implies
\begin{align*}
\mathcal{S}(z)\big(
   \TT_{k}^+ (u)\ot \TT_{m}^+ (v) \big)
=&\,\tr_{1,\ldots, m+k}\ts
S_{[k]}^1  S_{[m]}^2
A\cdotlr
\Big(B\ts
T_{[k]}^{+13}(u_{[k]})\ts
 C \ts
T_{[m]}^{+24}(v_{[m]}) \ts
B\ts(\vac\otimes \vac) \hspace{-2pt}\Big).
\end{align*}
Note that $S_{[m]}^2    S_{[k]}^1= S_{[k]}^1 S_{[m]}^2$, so we can use \eqref{zabc} to move $S_{[k]}^1$  to the right, thus getting
\begin{align*}
& \tr_{1,\ldots, m+k}\ts
  S_{[m]}^2
\wht{A} \cdotlr
\Big(\wht{B} \ts
S_{[k]}^1\ts
T_{[k]}^{+13}(u_{[k]})\ts
 C \ts
T_{[m]}^{+24}(v_{[m]}) \ts
B\ts(\vac\otimes \vac) \hspace{-2pt}\Big).
\end{align*}
Since
   $S_{[k]}^1 = (S_{[k]}^1)^2 $, this equals
\begin{align*}
\tr_{1,\ldots, m+k}\ts
  S_{[m]}^2
\wht{A} \cdotlr
\Big(\wht{B} \ts
(S_{[k]}^1)^2\ts
T_{[k]}^{+13}(u_{[k]})\ts
 C \ts
T_{[m]}^{+24}(v_{[m]}) \ts
B\ts(\vac\otimes \vac) \hspace{-2pt}\Big).
\end{align*}
Using \eqref{fusion2} and \eqref{zabc} we move one copy of $S_{[k]}^1$ to the left and another copy to the right:
\begin{align*}
\tr_{1,\ldots, m+k}\ts
 S_{[k]}^1 S_{[m]}^2 \ts
A \cdotlr
\Big(B \ts
\cev{T}_{[k]}^{+13}(u_{[k]})\ts
 \wht{C} \ts
T_{[m]}^{+24}(v_{[m]}) \ts
\wht{B}\ts S_{[k]}^1\ts (\vac\otimes \vac) \hspace{-2pt}\Big).
\end{align*}
By employing the cyclic property of the trace and then $(S_{[k]}^1)^2 =S_{[k]}^1$  we get
\begin{align*}
& \tr_{1,\ldots, m+k}\ts
 S_{[m]}^2 \ts
A \cdotlr
\Big(B \ts
\cev{T}_{[k]}^{+13}(u_{[k]})\ts
 \wht{C} \ts
T_{[m]}^{+24}(v_{[m]}) \ts
\wht{B}\ts (S_{[k]}^1)^2\ts (\vac\otimes \vac) \hspace{-2pt}\Big)\\
=\,& \tr_{1,\ldots, m+k}\ts
 S_{[m]}^2 \ts
A \cdotlr
\Big(B \ts
\cev{T}_{[k]}^{+13}(u_{[k]})\ts
 \wht{C} \ts
T_{[m]}^{+24}(v_{[m]}) \ts
\wht{B}\ts S_{[k]}^1\ts (\vac\otimes \vac) \hspace{-2pt}\Big).
\end{align*}
Next, we use equalities \eqref{fusion2} and \eqref{zabc} to move $S_{[k]}^1$ to the left:
\begin{align*}
\tr_{1,\ldots, m+k}\ts
 S_{[m]}^2 \ts
A \cdotlr
\Big(B \ts  S_{[k]}^1\ts
T_{[k]}^{+13}(u_{[k]})\ts
 C \ts
T_{[m]}^{+24}(v_{[m]}) \ts
B\ts (\vac\otimes \vac) \hspace{-2pt}\Big).
\end{align*}
Note that the tensor copies $1,\ldots ,k$ of $A$ commute with $S_{[m]}^2$, so this is equal to
$$
\tr_{1,\ldots, m+k}\ts
A \cdotlr
\Big( S_{[m]}^2 \ts B \ts  S_{[k]}^1\ts
T_{[k]}^{+13}(u_{[k]})\ts
 C \ts
T_{[m]}^{+24}(v_{[m]}) \ts
B\ts (\vac\otimes \vac) \hspace{-2pt}\Big).
$$
By using the cyclic property of the trace and then the first equality in \eqref{csyms} we get
\begin{align}
 & \tr_{1,\ldots, m+k}\ts
  S_{[m]}^2 \ts B \ts  S_{[k]}^1\ts
T_{[k]}^{+13}(u_{[k]})\ts
 C \ts
T_{[m]}^{+24}(v_{[m]}) \ts
B\ts A\ts (\vac\otimes \vac)\non  \\
=\,& \tr_{1,\ldots, m+k}\ts
  S_{[m]}^2 \ts B \ts  S_{[k]}^1\ts
T_{[k]}^{+13}(u_{[k]})\ts
 C \ts
T_{[m]}^{+24}(v_{[m]}) \ts (\vac\otimes \vac).  \non
\end{align}
Again, by the cyclic property of the trace,   this equals
\begin{align*}
& \tr_{1,\ldots, m+k}\ts
  S_{[m]}^2  \ts  S_{[k]}^1\ts
T_{[k]}^{+13}(u_{[k]})
 \left(B\cdotrl C\right)
T_{[m]}^{+24}(v_{[m]}) \ts (\vac\otimes \vac).
\end{align*}
Finally, as $S_{[m]}^2    S_{[k]}^1= S_{[k]}^1 S_{[m]}^2$, by employing the second equality in \eqref{csyms} we obtain
\begin{align*}
& \tr_{1,\ldots, m+k}\ts
 S_{[k]}^1\ts
T_{[k]}^{+13}(u_{[k]})
  S_{[m]}^2  \ts
T_{[m]}^{+24}(v_{[m]}) \ts (\vac\otimes \vac)
=
\TT_{k}^+ (u)\ot \TT_{m}^+ (v),
\end{align*}
as required.
\end{prf}

Let us  consider the classical limit of the series $\TT_{m}^+ (u)$. Recall \eqref{grading1} and define $\deg u=1$ and $\deg\partial_u =-1$.
Clearly, this extends \eqref{filtration} to the ascending filtration of the algebra $\Y^+(\g_N)[[u,\partial_u]]_{\text{fin}}$ which consists of all finite degree elements in $\Y^+(\g_N)[[u,\partial_u]]$. By Proposition \ref{pbwthm}, the corresponding graded algebra is $\left( \U(t^{-1}\g_N[t^{-1}])\ot_{\CC} \CC[[h]] \right) \hspace{-1pt}[[u,\partial_u]]_{\text{fin}}$ which  consists of all finite degree elements in $\left( \U(t^{-1}\g_N[t^{-1}])\ot_{\CC} \CC[[h]] \right) \hspace{-1pt}[[u,\partial_u]]$.
 Consider the element
\beq\label{cllim}
\tr_{1,\ldots ,m}\, S_{[m]}\left(1-T_1^+ (u)e^{- h\partial_u}\right)\left(1-T_2^+ (u)e^{- h\partial_u}\right)\ldots \left(1-T_m^+ (u)e^{- h\partial_u}\right) .
\eeq
Its degree equals $-m$, so it belongs  to the algebra $\Y^+(\g_N)[[u,\partial_u]]_{\text{fin}}$. By Proposition \ref{pbwthm}, its image in the corresponding graded algebra   equals
\beq\label{cllim4}
h^m \ts \tr_{1,\ldots ,m}\, S_{[m]}\left(\partial_u + F_1^+ (u) \right)\left(\partial_u +F_2^+ (u)\right)\ldots \left(\partial_u + F_m^+ (u) \right).
\eeq
On the other hand, by arguing as in \cite{JKMY}, we can express \eqref{cllim} as follows. First, using $e^{- h\partial_u} u^r =(u- h)^r e^{- h\partial_u} $ with $r\in\mathbb{Z}$  we observe that the given element is equal to
\beq\label{cllim2}
\tr_{1,\ldots ,m}\, S_{[m]}
\sum_{k=0}^m \ts
\sum_{1\leqslant i_1<\ldots < i_k\leqslant m}
(-1)^k\ts
T_{i_1}^+ (u)\ts
T_{i_2}^+ (u- h)\ldots
T_{i_k}^+ (u- (k-1)h)\ts
e^{- kh\partial_{u}}.
\eeq
Moreover, due to \cite[Lemma 4.1]{M}, we have
$$
\tr_m\, S_{[m]} = a_m\ts S_{[m-1]}\quad\text{for}\quad a_m=\pm \frac{(\pm N+m-3)(\pm N+2m-2)}{m(\pm N+2m -4)},
$$
where  the upper sign  corresponds to the orthogonal and the lower sign corresponds to the symplectic case. Conjugating the summands under the trace by  suitable elements of the symmetric group $\mathfrak{S}_m$ and using the first equality in \eqref{important} and the cyclic property of the trace we  write  \eqref{cllim2}  as
\begin{align*}
&\sum_{k=0}^m \ts
b_k\ts
\tr_{1,\ldots ,k}\, S_{[k]}\ts
T_{1}^+ (u)\ts
T_{2}^+ (u- h)\ldots
T_{k}^+ (u- (k-1)h)\ts
e^{- kh\partial_{u}}=\sum_{k=0}^m \ts
b_k\ts\TT_k^+ (u)\ts
e^{- kh\partial_{u}},
\end{align*}
where $b_k = (-1)^k a_{k+1}\ldots a_m  \binom{m}{k}$.
In the above equality we used alternative expression \eqref{intheform} for $\TT_k^+ (u) $ in the orthogonal case, so that our calculation treats both cases simultaneously.
Introduce the series
$$\Phi_{m}(u)=\sum_{r\geqslant 0}\Phi_{mr} u^r = h^{-m} \sum_{k=0}^m \ts b_k\ts\TT_k^+ (  u) \in \Vcr(\g_N)[[u]].$$
 Theorem \ref{thmmain4} implies
\begin{kor}\label{lagana}
All coefficients of $\Phi_{m}(u)$ belong to the center of the $h$-adic quantum vertex algebra $\Vcr(\g_N)$.
\end{kor}

Let $D \colon V_{\text{crit}}(\g_N)\to  V_{\text{crit}}(\g_N)$ be the translation operator defined on the universal affine vertex algebra $V_{\text{crit}}(\g_N)=V_{c_\text{crit}}(\g_N)$ at the critical level $c_{\text{crit}}$ by
$$D \vac=0\Fand [D, f_{ij} (-r)] = r f_{ij}(-r-1)\quad\text{for all}\quad i,j=1,\ldots N,\, r\geqslant 1,$$
where $\vac$ is the vacuum vector in $V_{\text{crit}}(\g_N)$.
By combining explicit formula \eqref{Doperator} for the operator
$\Dc\colon \Vc_{\text{crit}}(\g_N)\to  \Vc_{\text{crit}}(\g_N)$
 and Proposition \ref{pbwthm} one can easily verify that  map \eqref{pbw} acts as follows
\beq\label{observation9}
\Vc_{\text{crit}}(\g_N)\ni\Dc\ts  t_{i_1 j_1}^{(-r_1)}\ldots t_{i_s j_s}^{(-r_s)}\vac \,
\mapsto\,
D f_{i_1 j_1}(-r_1)\ldots f_{i_s j_s}(-r_s)\vac \in V_{\text{crit}}(\g_N)
\eeq
 for all $i_1,\ldots ,j_s=1,\ldots ,N$, $r_1,\ldots ,r_s\geqslant 1$ and $s\geqslant 0$, thus justifying our notation.

Denote by $\z(V_{\text{crit}}(\g_N))$ the  Feigin--Frenkel center \cite{FF}, i.e. the center of the  vertex algebra $V_{\text{crit}}(\g_N)$.
The complete sets $\left\{\phi_{2\ts 2},\phi_{44},\ldots,\phi_{2n\ts 2n}\right\}\subset\z(V_{\text{crit}}(\g_N))$ of Segal--Sugawara vectors for  $V_{\text{crit}}(\g_N)$, where $\g_N =\on_{2n+1}, \spn_{2n},\on_{2n}$, i.e. the sets such that
\beq\label{by14}
\z(V_{\text{crit}}(\g_N))=\CC[D^k \ts \phi_{2i\ts 2i}: i=1,\ldots ,n,\, k\geqslant 0]
\eeq
were constructed by Molev in \cite{M}. Extend the affine Lie algebra $\wht{\g}_N$ with the element $\tau$ so that the following commutation relations hold on $\wht{\g}_N\oplus\CC\tau$:
$$[\tau,C]=0\Fand [\tau, f_{ij} (-r)] = r f_{ij}(-r-1)\quad\text{for all}\quad i,j=1,\ldots N,\, r\in\ZZ. $$
The Segal--Sugawara vectors $\phi_{2i\ts 2i}\in V_{\text{crit}}(\g_N)\equiv \U(t^{-1}\g_N[t^{-1}])$ for  $i=1,2,\ldots, p_{\g_N}$, where
$$p_{\g_N}=\begin{cases}
n, &\text{if }\g_N=\on_{2n+1},\\
 \left\lfloor{n/2}\right \rfloor ,&\text{if }\g_N=\spn_{2n},\\
n-1 ,&\text{if }\g_N=\on_{2n},
\end{cases}
$$
 are found in \cite{M} as constant terms of the  polynomials
\beq\label{polys}
\tr_{1,\ldots ,m} \ts S_{[m]}
\left(\tau +F^+ (-1)_1\right)
\ldots
\left(\tau +F^+ (-1)_m\right)
=
\phi_{m\ts 0}\tau^m+\phi_{m\ts 1}\tau^{m-1}
+\ldots +\phi_{m\ts m}
\eeq
with $m=2,4,\ldots ,2p_{\g_N}$.
We now recover these vectors  by taking the classical limit of the constant terms of the series $\Phi_m (u)$.

\begin{pro}\label{algind}
The images  of the elements  $\Phi_{2i\ts 0}\in \z(\Vcr(\g_N))$   with respect to map \eqref{pbw} coincide  with the Segal--Sugawara vectors  $\phi_{2i\ts 2i}\in \z(V_{\text{crit}}(\g_N))$ for all $i=1,\ldots, p_{\g_N}$.
In particular, the images of the elements
$ \Phi_{2\ts 0} , \Phi_{4\ts 0} , \ldots , \Phi_{2n\ts 0} \in  \z(\Vcr(\on_{2n+1}))$
form a complete set of Segal--Sugawara vectors for the universal affine vertex algebra $V_{\text{crit}}(\on_{2n+1})$.
\end{pro}

\begin{prf}
Recall  \eqref{cllim4}.
The image of  the element $\Phi_{m\ts 0}\in \z(\Vcr(\g_N))$ with respect to map \eqref{pbw} is   found
 by moving all $\partial_u$ in
$$
\tr_{1,\ldots ,m}\, S_{[m]}\left(\partial_u + F_1^+ (u) \right)\left(\partial_u + F_2^+ (u) \right)\ldots \left(\partial_u + F_m^+ (u) \right)
$$
  to the right and then taking the constant term with respect to $u$ and $\partial_u$.
However, since
$$[\partial_u, F^+(u)]=\sum_{r \geqslant 1} r F(-r-1) u^{r-1}=\sum_{r \geqslant 1} [\tau,F(-r)] u^{r-1}=[\tau, F^+(u)],$$
it is clear from \eqref{polys} that this image   coincides with $\phi_{m\ts m}$, as required.
\end{prf}

Let $V$ be an  $h$-adic quantum vertex algebra with vacuum vector $\vac$. The product
\beq\label{centerprods}
a\cdot b =a_{-1} b\quad\text{for all }a,b\in \z(V)
\eeq
defines the structure of an  associative algebra with unit $\vac$
on $\z(V)$. Moreover, this algebra is equipped with a derivation defined as the restriction
of the   operator $\Dc$; see  \cite[Proposition 3.7]{JKMY}.
We   now    explicitly describe the center of the $h$-adic quantum vertex algebra $ \Vcr(\on_{2n+1})$.

\begin{thm}\label{FFthm}
The algebra $\z(\Vcr(\on_{2n+1}))$  coincides with the
 $h$-adically completed polynomial algebra in infinitely many variables,
\beq\label{polyalg}\z(\Vcr(\on_{2n+1}))=\CC[\Dc^k\ts \Phi_{2i\ts 0} : i=1,\ldots ,n,\, k\geqslant 0][[h]].
\eeq
\end{thm}

\begin{prf}
First, we note that by Corollary \ref{lagana} and \cite[Proposition 3.7]{JKMY} all elements $\Dc^k \Phi_{2i\ts 0}$ with $k\geqslant 0$ belong to the center.
Furthermore, all these elements are algebraically independent.
Indeed, by Proposition \ref{algind} and \eqref{observation9}
 their images with respect to map \eqref{pbw} coincide with elements of $V_{\text{crit}}(\on_{2n+1})$, which are algebraically independent
 due to \cite{M}. Denote by $\Pc$ the $h$-adically completed subalgebra of $\z(\Vcr(\on_{2n+1})) $ generated by all $\Dc^k \Phi_{2i\ts 0}$. By  \cite[Proposition 3.7]{JKMY} we have $\Pc\subseteq \z(\Vcr(\on_{2n+1}))$.

The $\CC[[h]]$-module $\Vcr(\on_{2n+1})$ is topologically free, so it can be written as $\Vc_0[[h]]$ for some complex vector space $\Vc_0$.
Let $v$ be an arbitrary element of the center $\z(\Vcr(\on_{2n+1}))$. We will prove by induction that for every nonnegative integer  $r$ there exists an element
$p\in\Pc$
such that
$v-p$
belongs to $h^r \Vcr(\on_{2n+1})$.

First, note that the aforementioned statement  holds trivially for $r=0$.
Suppose that for some integer $r\geqslant 0$ we have
$$v-p =h^r v_r +h^{r+1} v_{r+1}+\ldots \qquad\text{for some }p\in\Pc\text{ and }v_s\in\Vc_0 .$$
Clearly, both $v-p$ and  $v_r+h v_{r+1}+h^2 v_{r+2}+\ldots$ belong to the center $\z(\Vcr(\on_{2n+1}))$.
Write $v_r$ as the sum $v_r=v_{r,1}+\ldots + v_{r,l_r}$, so that all  $v_{r,j}$ are homogeneous with respect to degree operator \eqref{grading1}. Set $d_r=\min\textstyle\left\{\deg v_{r,j}:j=1,\ldots ,l_r\right\}$.

Choose  any homogeneous element
$$w =h^{s_1}v_{s_1} +h^{s_2}v_{s_2}+\ldots\qquad\text{with }s_{j+1}>s_j\geqslant r\text{ for all }j$$
such that $\deg(v-p-w)< \deg(v-p)=\deg w=\deg v_{s_j}$ for all $j$.
By combining \eqref{centar} and Proposition \ref{pbwthm} we conclude that the images      $\bar{v}_{ s_j}$ of  $v_{ s_j}$ with respect to map \eqref{pbw} belong to the center  of the vertex algebra $V_{\text{crit}}(\on_{2n+1})$ for all $j\geqslant 0$.
Hence, due to \eqref{by14}, there exist  polynomials $\bar{p}_j$ in the variables $D^k  \phi_{2i\ts 2i}$ such that $\bar{v}_{ s_j}=\bar{p}_j$ for all $j $.
Let $p_j$ be the polynomials obtained from $\bar{p}_j$ by replacing the variables $D^k  \phi_{2i\ts 2i}$ with the respective variables $\Dc^k\ts \Phi_{2i\ts 0}$.
The element $v-p-\sum_{j} h^{s_j}p_j$ belongs to the center  $\z(\Vcr(\on_{2n+1}))$ and its degree is strictly less than $\deg(v-p)=\deg w$.

Let
\beq\label{startingwith}
v-p-\sum_{j} h^{s_j}p_j=h^r v_r' +h^{r+1} v_{r+1}'+\ldots \qquad\text{for some } v_s'\in\Vc_0 .
\eeq
Write $v_r'$ as a sum $v'_r=v'_{r,1}+\ldots + v'_{r,l'_r}$, so that all elements $v'_{r,j} $ are homogeneous with respect to degree operator \eqref{grading1}. Set $d_r'=\min\textstyle\left\{\deg v'_{r,j}:j=1,\ldots ,l_r'\right\}$. Observe that $d_r'\geqslant d_r$ because lower degree terms, with respect to  \eqref{grading1}, in  all elements $\Dc^k \Phi_{2i\ts 0}$ come up multiplied by a positive power of $h$.
We can continue to  repeat the same procedure, now starting with element \eqref{startingwith}, for an appropriate number of times. As we demonstrated, in each  step   the degree of the left hand side is reduced  while  the degree of the lowest degree term of the coefficient of $h^r$ on the right hand side does not decrease. Therefore, after finitely many steps we  end up with the expression of the form
$$
v-p-\sum_{j} h^{t_j}q_j=h^{r+1} v_{r+1}''+h^{r+2} v_{r+2}''+\ldots \qquad\text{for some } v_s''\in\Vc_0\text{ and }q_j\in\Pc,
$$
 thus finishing the inductive step. Note that the sequence $(t_j)$ is strictly increasing, so that the sum $\sum_{j} h^{t_j}q_j$ is indeed a well-defined element of $\Pc$.  Therefore, we proved   that for any  $v\in\z(\Vcr(\on_{2n+1}))$ and for any integer $r\geqslant 0$ there exists
$p\in\Pc$
such that
$v-p$
belongs to $h^r \Vcr(\on_{2n+1})$. This implies $\Pc = \Vcr(\on_{2n+1})$.

It remains to prove that the algebra $\Pc$  is commutative.
By \cite[Proposition 3.8]{JKMY}  the center of every $h$-adic quantum vertex algebra is $\Sc$-commutative, i.e. we have
$$Y(z_1) (1\ot Y(z_2))\Sc(z_1 -z_2) (a\ot b)
=Y(b,z_2)Y(a,z_1)\quad\text{for all }a,b\in\z(\Vcr(\on_{2n+1})).
$$
By setting $a=\Phi_{2i\ts 0}$, $b=\Phi_{2j\ts 0}$ and using Theorem \ref{thmmain4} we find
\beq\label{derivatives4}
Y(\Phi_{2i\ts 0},z_1) \ts Y(\Phi_{2j\ts 0},z_2)
=Y(\Phi_{2j\ts 0},z_2)\ts Y(\Phi_{2i\ts 0},z_1).
\eeq
Due to \cite[Lemma 2.13]{Li} we have
$$Y(\Dc v,z)=\frac{d}{dz} Y(v,z)\quad\text{for all } v\in \Vcr(\on_{2n+1}),$$
so by  applying the partial derivatives $\partial^k/\partial z_1^k$ and $\partial^m/\partial z_1^m$ to \eqref{derivatives4} we get
$$
Y(\Dc^k\ts\Phi_{2i\ts 0},z_1) \ts Y(\Dc^m\ts\Phi_{2j\ts 0},z_2)
=Y(\Dc^m\ts\Phi_{2j\ts 0},z_2)\ts Y(\Dc^k\ts\Phi_{2i\ts 0},z_1).
$$
By arguing as in the proof of \cite[Proposition 3.4]{K9} one can  prove that this implies
$$
Y(a,z_1) \ts Y(b,z_2)
=Y(b,z_2)\ts Y(a,z_1)\quad\text{for all }  a,b\in\z(\Vcr(\on_{2n+1})).
$$
Finally,  applying this equality to the vacuum vector $\vac$ and then taking the constant terms with respect to the variables $z_1$ and $z_2$ we find
$a\cdot b=b\cdot a$, as required.
\end{prf}

It is worth to single out the commutativity property of the restriction of  vertex operator map \eqref{qva1} on the center, which was obtained in the proof of Theorem \ref{FFthm}.

\begin{kor}\label{korsat}
For any $a,b\in\z(\Vcr(\on_{2n+1}))$ we have
$$
Y(a,z_1) \ts Y(b,z_2)
=Y(b,z_2)\ts Y(a,z_1).
$$
\end{kor}
By arguing as in the proof of Theorem \ref{FFthm}, one  obtains the following partial result on the quantum center in types $C$ and $D$.
\begin{kor}\label{korsat2}
The algebra $\z(\Vcr(\g_N))$ with $\g_N=\spn_{2n},\on_{2n}$  contains the
 $h$-adically completed polynomial algebra in infinitely many variables
\beq\label{CD}
\CC[\Dc^k\ts \Phi_{2i\ts 0} : i=1,\ldots ,p_{\g_N},\, k\geqslant 0][[h]].
\eeq
\end{kor}

In particular, the commutativity of the restriction of the vertex operator map $Y(z)$ on \eqref{CD} is established by arguing as in the last part of the proof of Theorem \ref{FFthm}.
Recall that by \eqref{centar} the center $\z(\Vcr(\g_N))$ consists of $\Y(\g_N)$-invariants, so that  we have
$$
Y\big(T_{[n]}^+ (u)\vac,z\big) v=T_{[n]}^+ (u|z)\ts v\quad\text{for all }v\in \z(\Vcr(\g_N)).
$$
By combining these  two observations with Corollaries \ref{korsat} and \ref{korsat2} we obtain commutative subalgebras of the dual Yangian in types $B$, $C$ and $D$, as suggested by \cite[Remark 11.2.5]{Mnew}.
\begin{kor}
The coefficients of the series $\TT_m^+(u)$, where $m= 1,\ldots ,N$ in the orthogonal case and $m = 1,...,N/2$ in the
symplectic case, generate a commutative subalgebra of the dual Yangian $\Y^+(\g_N)$.
\end{kor}

\section{Central elements of the completed double Yangian  at the critical level}

We now employ \eqref{tcenter} to obtain explicit formulae for families of central elements of the appropriately completed double Yangian  at the critical level.
Introduce the completion of the double Yangian $\DY_c(\g_N)$ at the level $c\in\CC$ as the inverse limit
$$
\DYtld_c(\g_N)=\lim_{\longleftarrow} \,\DY_c(\g_N)/\, \I_p,
$$
where $p\geqslant 1$ and $\I_p$ denotes the $h$-adically completed left ideal of $\DY_c(\g_N)$ generated by all elements $t_{ij}^{(r)}$ with $r\geqslant p$.

From now on, we consider the completed double Yangian $\DYtldc(\g_N)$ at the critical level $c_\textrm{crit}=-2\kp/\sigma$. Introduce the series in
$\DYtldc(\g_N)[[u^{\pm 1}]]$ by
\begin{align}\label{oznaka1}
\TT_m(u)=&\,\tr_{1,\ldots ,m} \ts S_{[m]} T_{[m]}^+ (u_{[m]})T_{[m]} (u_{[m]}-h\kp)^{-1}\non\\
=&\,\tr_{1,\ldots ,m} \ts S_{[m]} T_{1}^+ (\wvr{u}_1)\ldots T_{m}^+ (\wvr{u}_m)
T_{m} (\wvr{u}_m-h\kp)^{-1}\ldots
T_{1} (\wvr{u}_1-h\kp)^{-1},
\end{align}
where $m=1,\ldots , N$ in the orthogonal case,  $m=1,\ldots, N/2$ in the symplectic case, the trace is taken over all $m$ copies of $\ndo\CC^N$ and
$u_{[m]}=(\wvr{u}_{1},\ldots ,\wvr{u}_m)$
is given by \eqref{uem1}. The proof of the next theorem is analogous to the proof of \cite[Theorem 3.2]{FJMR}.

\begin{thm}\label{lastthm}
All coefficients of $\TT_m(u)$ belong to the center of the completed double Yangian $\DYtldc(\g_N)$.
\end{thm}

\begin{prf}
It is sufficient to verify the equalities
\beq\label{finalXXX}
T  (v_0)\TT_m(u)=\TT_m(u)T  (v_0)\fand T^+ (v_{m+1})\TT_m(u)=\TT_m(u)T^+ (v_{m+1}).
\eeq

Let us prove the first equality in \eqref{finalXXX}.
Label the tensor copies as follows:
$$
\overbrace{\ndo\CC^N}^{0} \ot \overbrace{(\ndo \CC^N)^{\ot m}}^{1} \ot \overbrace{\DYtldc(\g_N)}^{2}.
$$
The elements
$$A=\R_{1m}^{01}(v_0-h\kp|u_{[m]})^{-1}\fand
B=\R_{1m}^{01}(v_0+h\kp|u_{[m]})$$
satisfy
\beq\label{useful1}
S_{[m]}^1 A=\cev{A}S_{[m]}^1,\qquad S_{[m]}^1 B=\cev{B}S_{[m]}^1\Fand A\cdotlr B=1.
\eeq
 Indeed, the first two equalities follow from \eqref{fusion3} while the third equality is a consequence of crossing symmetry properties \eqref{csym} and \eqref{csym2}.

By applying $T  (v_0)$ on $\TT_m(u)$ and using defining relations \eqref{DX1} and \eqref{DX3} we find
\begin{align*}
T (v_0)\TT_m(u)=&\,
T_{[1]}^0 (v_0) \,\tr_{1,\ldots ,m} \ts S_{[m]}^1\ts T_{[m]}^{+12} (u_{[m]})\ts T_{[m]}^{12} (u_{[m]}-h\kp)^{-1} \\
=&\, \tr_{1,\ldots ,m} \ts S_{[m]}^1\ts T_{[1]}^{02} (v_0) \ts T_{[m]}^{+12} (u_{[m]})\ts T_{[m]}^{12} (u_{[m]}-h\kp)^{-1}\\
=&\, \tr_{1,\ldots ,m} \ts S_{[m]}^1 A\ts T_{[m]}^{+12} (u_{[m]})\ts T_{[1]}^{02} (v_0) \ts B\ts T_{[m]}^{12} (u_{[m]}-h\kp)^{-1}\\
=&\, \tr_{1,\ldots ,m} \ts S_{[m]}^1 A\ts T_{[m]}^{+12} (u_{[m]})\ts   T_{[m]}^{12} (u_{[m]}-h\kp)^{-1}\ts B \ts T_{[1]}^{02} (v_0).
\end{align*}
Therefore, by employing the first equality in \eqref{useful1} we obtain
\beq\label{rhs1}
T  (v_0)\TT_m(u)=\tr_{1,\ldots ,m} \ts \cev{A}\ts S_{[m]}^1  T_{[m]}^{+12} (u_{[m]})\ts   T_{[m]}^{12} (u_{[m]}-h\kp)^{-1}  B \ts\ts T^{02}_{[1]}  (v_0).
\eeq

Since $S_{[m]}=(S_{[m]})^2$ we can write
$$\tr_{1,\ldots ,m} \ts \cev{A}\ts S_{[m]}^1  T_{[m]}^{+12} (u_{[m]})\ts   T_{[m]}^{12} (u_{[m]}-h\kp)^{-1}  B $$
as
$$\tr_{1,\ldots ,m} \ts \cev{A}\ts (S_{[m]}^1)^2  \ts T_{[m]}^{+12} (u_{[m]})\ts   T_{[m]}^{12} (u_{[m]}-h\kp)^{-1}  B. $$
We now use  \eqref{fusion22}, \eqref{fusion2} and the first two equalities in \eqref{useful1} to move one copy of the symmetrizer $S_{[m]}$ to the left and another copy to the right, thus getting
$$\tr_{1,\ldots ,m} \ts S_{[m]}^1\ts A\ts  \cev{T}_{[m]}^{+12} (u_{[m]})\ts   \cev{T}_{[m]}^{12} (u_{[m]}-h\kp)^{-1}  \cev{B}\ts S_{[m]}^1. $$
By the cyclic property of the trace and $(S_{[m]})^2=S_{[m]}$ this equals
$$\tr_{1,\ldots ,m} \ts A\ts  \cev{T}_{[m]}^{+12} (u_{[m]})\ts   \cev{T}_{[m]}^{12} (u_{[m]}-h\kp)^{-1}  \cev{B}\ts S_{[m]}^1. $$
Next, we employ \eqref{fusion22}, \eqref{fusion2} and the second equality in \eqref{useful1} to move the symmetrizer to the left, thus getting
$$\tr_{1,\ldots ,m} \ts A\ts S_{[m]}^1\ts  T_{[m]}^{+12} (u_{[m]})\ts   T_{[m]}^{12} (u_{[m]}-h\kp)^{-1}  B . $$
Finally, by the cyclic property of the trace and the third equality in \eqref{useful1} this is equal to
$$\tr_{1,\ldots ,m} \ts S_{[m]}^1\ts  T_{[m]}^{+12} (u_{[m]})\ts   T_{[m]}^{12} (u_{[m]}-h\kp)^{-1}=\TT_m(u)   . $$
Therefore, we conclude that the right hand side in \eqref{rhs1} equals  $\TT_m(u)T  (v_0)$, as required.

Let us prove the second equality in \eqref{finalXXX}.
Label the tensor copies as follows:
$$
 \overbrace{(\ndo \CC^N)^{\ot m}}^{1} \ot\overbrace{\ndo\CC^N}^{2} \ot \overbrace{\DYtldc(\g_N)}^{3}.
$$
The elements
$$X=\R_{m1}^{12}(u_{[m]}|v_{m+1})\fand
Z=\R_{m1}^{12}(u_{[m]}|v_{m+1}+2h\kp)^{-1}.$$
satisfy
\beq\label{useful2}
S_{[m]}^1X=\cev{X}S_{[m]}^1,\qquad S_{[m]}^1Z=\cev{Z}S_{[m]}^1 \Fand X\cdotrl Z=1.
\eeq
As with \eqref{fusion3}, the first two equalities in \eqref{useful2} can be proved by using Yang--Baxter equation \eqref{YBE} and fusion procedure \eqref{fusion}.
The third equality is a consequence of crossing symmetry properties \eqref{csym} and \eqref{csym2}.

By applying $T^+ (v_{m+1})$ on $\TT_m(u)$ and using defining relations \eqref{DX2} and \eqref{DX3} we find
\begin{align*}
T^+ (v_{m+1})\TT_m(u)=&\,
T_{[1]}^{+23} (v_{m+1})\ts \tr_{1,\ldots ,m} \ts S_{[m]}^1\ts T_{[m]}^{+13} (u_{[m]})\ts T_{[m]}^{13} (u_{[m]}-h\kp)^{-1}\\
=&\, \tr_{1,\ldots ,m} \ts S_{[m]}^1\ts T_{[1]}^{+23} (v_{m+1}) \ts T_{[m]}^{+13} (u_{[m]})\ts T_{[m]}^{13} (u_{[m]}-h\kp)^{-1}\\
=&\, \tr_{1,\ldots ,m} \ts S_{[m]}^1 X\ts T_{[m]}^{+13} (u_{[m]})\ts T_{[1]}^{+23} (v_{m+1})\ts X^{-1}\ts T_{[m]}^{13} (u_{[m]}-h\kp)^{-1}\\
=&\, \tr_{1,\ldots ,m} \ts S_{[m]}^1 X\ts T_{[m]}^{+13} (u_{[m]})\ts   T_{[m]}^{13} (u_{[m]}-h\kp)^{-1}\ts Z \ts T_{[1]}^{+23} (v_{m+1}).
\end{align*}
Therefore, by employing the first equality in \eqref{useful2} we obtain
\beq\label{rhs2}
T^+ (v_{m+1})\TT_m(u)=\tr_{1,\ldots ,m} \ts \cev{X}\ts S_{[m]}^1  T_{[m]}^{+13} (u_{[m]})\ts   T_{[m]}^{13} (u_{[m]}-h\kp)^{-1}  Z \ts\ts T_{[1]}^{+23} (v_{m+1}).
\eeq

Since $S_{[m]}=(S_{[m]})^2$ we can write
$$\tr_{1,\ldots ,m} \ts \cev{X}\ts S_{[m]}^1  T_{[m]}^{+13} (u_{[m]})\ts   T_{[m]}^{13} (u_{[m]}-h\kp)^{-1}  Z $$
as
$$\tr_{1,\ldots ,m} \ts \cev{X}\ts (S_{[m]}^1)^2 \ts  T_{[m]}^{+13} (u_{[m]})\ts   T_{[m]}^{13} (u_{[m]}-h\kp)^{-1}  Z .$$
We now use  \eqref{fusion22}, \eqref{fusion2} and the first two equalities in \eqref{useful2} to move one copy of the symmetrizer $S_{[m]}$ to the left and another copy to the right, thus getting
$$\tr_{1,\ldots ,m} \ts S_{[m]}^1\ts X\ts  \cev{T}_{[m]}^{+13} (u_{[m]})\ts   \cev{T}_{[m]}^{13} (u_{[m]}-h\kp)^{-1}  \cev{Z}\ts S_{[m]}^1. $$
By the cyclic property of the trace and $(S_{[m]})^2=S_{[m]}$ this equals
$$\tr_{1,\ldots ,m} \ts X\ts  \cev{T}_{[m]}^{+13} (u_{[m]})\ts   \cev{T}_{[m]}^{13} (u_{[m]}-h\kp)^{-1}  \cev{Z}\ts S_{[m]}^1. $$
Next, we employ \eqref{fusion22}, \eqref{fusion2} and the second equality in \eqref{useful2} to move the symmetrizer to the left, thus getting
$$\tr_{1,\ldots ,m} \ts X\ts S_{[m]}^1\ts  T_{[m]}^{+13} (u_{[m]})\ts   T_{[m]}^{13} (u_{[m]}-h\kp)^{-1}  Z . $$
Finally, by the cyclic property of the trace and the third equality in \eqref{useful2} this is equal to
$$\tr_{1,\ldots ,m} \ts S_{[m]}^1\ts  T_{[m]}^{+13} (u_{[m]})\ts   T_{[m]}^{13} (u_{[m]}-h\kp)^{-1}=\TT_m(u)   . $$
This implies  that the right hand side in \eqref{rhs2} is equal to  $\TT_m(u)T^+ (v_{m+1})$, as required.
Therefore,   we proved both equalities in \eqref{finalXXX}, so the theorem   follows.
\end{prf}

In the end, it is worth noting the following equality
$$\TT_m(u)=Y(\TT_m^+(0)\vac,u)$$
for operators on $\Vcr(\g_N)$, which suggests the form of formulae \eqref{oznaka1}, as well as the close    connection between  double Yangians and the corresponding $h$-adic quantum vertex algebras which is yet to be investigated.


\section*{Acknowledgement}
The first author is partially supported by the QuantiXLie Centre of Excellence, a project cofinanced by the Croatian Government and European Union through the European Regional Development Fund - the Competitiveness and Cohesion Operational Programme  (Grant KK.01.1.1.01.0004). The second
author is supported by Simons Foundation grant no. 523868 and National Natural Science Foundation of China grant no. 11531004.

\end{document}